\documentclass[12pt]{article}

\usepackage{latexsym}
\usepackage{amsfonts}
\usepackage{amssymb}
\usepackage{amsmath}
\usepackage{mathrsfs}
\usepackage{bbm}
\usepackage{dsfont}
\usepackage{hyperref}

\newcommand{\be}{\begin{equation}}
\newcommand{\ee}{\end{equation}}
\newcommand{\bea}{\begin{eqnarray}}
\newcommand{\eea}{\end{eqnarray}}
\newcommand{\bean}{\begin{eqnarray*}}
\newcommand{\eean}{\end{eqnarray*}}
\newcommand{\brray}{\begin{array}}
\newcommand{\erray}{\end{array}}

\newtheorem{dfn}{Definition}[section]
\newtheorem{thm}[dfn]{Theorem}
\newtheorem{lmma}[dfn]{Lemma}
\newtheorem{ppsn}[dfn]{Proposition}
\newtheorem{crlre}[dfn]{Corollary}
\newtheorem{xmpl}[dfn]{Example}
\newtheorem{rmrk}[dfn]{Remark}

\newcommand{\bdfn}{\begin{dfn}\rm}
\newcommand{\bthm}{\begin{thm}}
\newcommand{\blmma}{\begin{lmma}}
\newcommand{\bppsn}{\begin{ppsn}}
\newcommand{\bcrlre}{\begin{crlre}}
\newcommand{\bxmpl}{\begin{xmpl}}
\newcommand{\brmrk}{\begin{rmrk}\rm}

\newcommand{\edfn}{\end{dfn}}
\newcommand{\ethm}{\end{thm}}
\newcommand{\elmma}{\end{lmma}}
\newcommand{\eppsn}{\end{ppsn}}
\newcommand{\ecrlre}{\end{crlre}}
\newcommand{\exmpl}{\end{xmpl}}
\newcommand{\ermrk}{\end{rmrk}}

\newcommand{\clh}{\mathcal{H}}

\newcommand{\clk}{\mathcal{K}}

\newcommand{\clg}{\mathcal{G}}

\author{Anbu Arjunan and S. Sundar}
\title{CCR flows associated to closed convex cones}

\begin{document}
\maketitle
\begin{abstract}
 Let $P$ be a closed convex cone in $\mathbb{R}^{d}$ which we assume to be spanning and pointed i.e. $P-P=\mathbb{R}^{d}$ and $P \cap -P=\{0\}$. In this article, we consider CCR flows over $P$ 
 associated to isometric representations that arise out of $P$-invariant closed subsets, also called as $P$-modules, of $\mathbb{R}^{d}$.  We show that for two $P$-modules the associated CCR
 flows are cocycle conjugate if and only if the modules are translates of each other. 
 \end{abstract}
\noindent {\bf AMS Classification No. :} {Primary 46L55; Secondary 46L99.}  

\tableofcontents

\section{Introduction}

The theory of $E_0$-semigroups initiated by Powers and further developed extensively by Arveson is approximately three decades old. 
We refer the reader to the beautiful monograph \cite{Arveson} for the history, the development and the literature on $E_0$-semigroups. 
In this long introduction, we explain the problem considered in this paper, collect the preliminaries required and explain the techniques behind
the proof of our main theorem. Let $\clh$ be an infinite dimensional  separable Hilbert space. We  denote the algebra of bounded operators on $\clh$ by $B(\clh)$. By an $E_0$-semigroup on $B(\clh)$, one means a $1$-parameter semigroup $\alpha:=\{\alpha_{t}\}_{t \geq 0}$ of unital normal $*$-endomorphisms of $B(\clh)$.  However nothing prevents us from considering semigroups of endomorphisms on $B(\clh)$ indexed by more general semigroups. 

The authors in collobaration with others (\cite{Murugan_Sundar_continuous}, \cite{Anbu}) have considered $E_0$-semigroups over closed convex cones. 
In this paper, we analyse the basic examples of Arveson's theory i.e. CCR flows associated to modules over cones. We hope that the reader will be convinced by the end of this paper that this is not merely for the sake of  generalisation and there are some interesting connections to groupoid $C^{*}$-algebras and in particular to the groupoid approach to topological semigroup $C^{*}$-algebras which was first systematically explored by Muhly and Renault  in \cite{Renault_Muhly}. 

We fix notation that will be used throughout this paper. The norm that we use on $\mathbb{R}^{d}$ is always the usual Euclidean norm.  Let $P \subset \mathbb{R}^{d}$ be a closed convex cone. We assume $P$ is pointed i.e. $P \cap -P=\{0\}$. By restricting ourselves to the vector space generated by $P$, there is no loss of generality in assuming that  $P$ is spanning i.e. $P-P=\mathbb{R}^{d}$. The interior of $P$ will be denoted by $\Omega$. Then $\Omega$ is dense in $P$.  For a proof of this, we refer the reader to Lemma 3.1 of \cite{Murugan_Sundar}. It is also clear that $\Omega$ spans $\mathbb{R}^{d}$. For $x,y \in \mathbb{R}^{d}$, we write $x\geq y$ if $x-y \in P$ and write $x>y$ if $x-y \in \Omega$. We have the following \textbf{Archimedean principle}: Given $x \in \mathbb{R}^{d}$ and $a \in \Omega$, there exists a positive integer $n_0$ such that $n_0a>x$. For a proof of this fact, we refer the reader to Lemma 3.1 of \cite{Murugan_Sundar_continuous}

Let us review the definitions of $E_0$-semigroups and some results from \cite{Anbu}. Let $\clh$ be an infinite dimensional separable Hilbert space. By an \emph{$E_0$-semigroup, over $P$}, on $B(\clh)$, we mean a family $\alpha:=\{\alpha_{x}\}_{x \in P}$ of normal unital $*$-endomorphisms of $B(\clh)$ such that $\alpha_{x} \circ \alpha_{y}=\alpha_{x+y}$ satisfying the following continuity condition:  For $A \in B(\clh)$ and $\xi,\eta \in \clh$, the map $P \ni x \to \langle \alpha_{x}(A)\xi|\eta \rangle$ is continuous. 

We consider two $E_0$-semigroups acting on different Hilbert spaces to be \emph{isomorphic}  if they are unitarily equivalent. The precise definition is as follows.
Let $\clk$ be an infinite dimensional separable Hilbert space and $U:\clh \to \clk$ be a unitary.  We denote the map $B(\clh) \ni T \to UTU^{*} \in B(\clk)$  by $Ad(U)$. Let $\alpha:=\{\alpha_{x}\}_{x \in P}$ and $\beta:=\{\beta_{x}\}_{x \in P}$ be $E_{0}$-semigroups acting on $B(\clh)$ and $B(\clk)$ respectively. We say that \emph{$\alpha$ is conjugate to $\beta$} if there exists a unitary $U:\clh \to \clk$ such that for every $x \in P$, $\beta_{x}=Ad(U)\circ \alpha_{x} \circ Ad(U)^{*}$. 

Let $\alpha:=\{\alpha_{x}\}_{x \in P}$ be an $E_0$-semigroup on $B(\clh)$. By an \emph{$\alpha$-cocycle}, we mean a strongly continuous family of unitaries $\{U_{x}\}_{x \in P}$ such that $U_{x}\alpha_{x}(U_y)=U_{x+y}$. If $U:=\{U_{x}\}_{x \in P}$ is an $\alpha$-cocycle, it is straightforward to check that $\{Ad(U_{x})\circ \alpha_{x}\}_{x \in P}$ is an $E_0$-semigroup. Such an $E_0$-semigroup is called a \emph{cocycle perturbation} of $\alpha$. Let $\beta$ be an $E_0$-semigroup acting on a possibly different Hilbert space say $\clk$. We say that $\beta$ is \emph{cocycle conjugate} to $\alpha$ if a conjugate of $\beta$ is a cocycle perturabation of $\alpha$. 

One natural operation that one can do with $E_0$-semigroups is the tensor product operation. Let $\alpha$ and $\beta$ be $E_0$-semigroups  on $B(\clh)$ and $B(\clk)$ respectively. Then there exists a unique $E_0$-semigroup, denoted $\alpha \otimes \beta$, on $B(\clh \otimes \clk)$ such that for $x \in P$, $A \in B(\clh)$ and $B \in B(\clk)$, \[(\alpha \otimes \beta)_{x}(A \otimes B)=\alpha_{x}(A) \otimes \beta_{x}(B).\]  For a proof of this fact, we refer the reader to the paragraph preceeding Remark 4.8 of \cite{Anbu}. It is  routine to verify that if $\beta$ is cocycle conjugate to $\gamma$ then $\alpha \otimes \beta $ is cocycle conjugate to $\alpha \otimes \gamma$. 

As with any mathematical structures, the first question is to know whether there are enough examples and, if possible, how to classfiy them. It is beyond the scope of the present paper to offer a complete classification of $E_0$-semigroups. This question is still open even in the $1$-dimensional case. We present here a class of examples that we call \emph{the CCR flows associated to $P$-modules} and we classify them completely.

First let us recall the notion of Weyl operators on the symmetric Fock space of $\clh$. Let $\Gamma(\clh)$ be the symmetric Fock space. For $u \in \clh$, let \[
e(u):=\sum_{n=0}^{\infty}\frac{u^{ \otimes n}}{\sqrt{n!}}.\] The set of vectors $\{e(u): u \in \clh\}$ is called the set of exponential vectors. We have the following. 
\begin{enumerate}
\item[(1)] For $u,v \in \clh$, $\langle e(u)|e(v) \rangle = e^{\langle u|v \rangle}$.
\item[(2)] The set $\{e(u): u \in \clh\}$ is total in $\Gamma(\clh)$.
\item[(3)] Any finite subset of $\{e(u):u \in \clh\}$ is linearly independent. 
\end{enumerate}
For $u \in \clh$, there exists a unique unitary, $W(u)$ on $\Gamma(\clh)$, whose action on the exponential vectors is given by the following formula:
\[
W(u)e(v):=e^{-\frac{||u||^{2}}{2}-\langle u|v \rangle}e(u+v). 
\]
The operators $\{W(u): u \in \clh\}$ are called \emph{the Weyl operators}. The Weyl operators satisfy the following canonical commutation relation. For $u,v \in \clh$, 
\[
W(u)W(v)=e^{-i Im  \langle u|v \rangle}W(u+v)
\] 
where $Im\langle u|v \rangle$ denotes the imaginary part of $\langle u|v \rangle$.  The linear span of the Weyl operators $\{W(u): u \in \clh\}$ forms a unital $*$-subalgebra of $B(\Gamma(\clh))$ whose strong closure is 
$B(\Gamma(\clh))$. 

For a unitary $U$ on $\clh$, there exists a unique unitary, denoted $\Gamma(U)$, on $\Gamma(\clh)$, whose action on the exponential vectors is given by 
\[
\Gamma(U)e(v):=e(Uv).
\]
The unitary $\Gamma(U)$ is called \emph{the second quantisation of $U$}. For a unitary $U$ on $\clh$ and $u \in \clh$, we have the relation $\Gamma(U)W(u)\Gamma(U)^{-1}=W(Uu)$. 

Let $\clh_1$ and $\clh_2$ be Hilbert spaces. The map \[\Gamma(\clh_1 \oplus \clh_2) \ni e(u_1 \oplus u_2) \to e(u_1) \otimes e(u_2) \in \Gamma(\clh_1) \otimes \Gamma(\clh_2)\]extends to a unitary. Via this unitary, we always identify $\Gamma(\clh_1 \oplus \clh_2)$ with $\Gamma(\clh_1) \otimes \Gamma( \clh_2)$. Under this identification, we have the equality $W(u_1 \oplus u_2)=W(u_1)\otimes W(u_2)$ for $u_1 \in \clh_1$ and $u_2 \in \clh_2$. 
For proofs of all the assertions made so far, we refer the reader to the book \cite{KRP}.

\begin{dfn}
By a strongly continuous isometric representation of $P$ on $\clh$, we mean a  map $V:P \to B(\clh)$ such that 
\begin{enumerate}
\item[(1)] for $x,y \in P$, $V_{x}V_{y}=V_{x+y}$,
\item[(2)] for $x \in P$, $V_{x}$ is an isometry, and
\item[(3)] for $\xi \in \clh$, the map $P \ni x \to V_{x}\xi \in \clh$ is continuous. 
\end{enumerate}
\end{dfn}

Let $V:P \to B(\clh)$ be a strongly continuous isometric representation. Then there exists a unique $E_0$-semigroup on $B(\Gamma(\clh))$ denoted $\alpha^{V}:=\{\alpha_{x}\}_{x \in P}$ such that for $x \in P$ and $u \in \clh$, \[\alpha_{x}(W(u))=W(V_{x}u).\]
For the existence of the $E_0$-semigroup $\alpha^{V}$, we refer the reader to Prop. 4.7 of \cite{Anbu}. We call $\alpha^{V}$ \emph{the CCR flow associated to the isometric representation $V$}. The association $V \to \alpha^{V}$ converts the direct sum of isometric representations to tensor product of $E_0$-semigroups. That is, let $V_1$ and $V_2$ be isometric representations of $P$ on  Hilbert spaces $\clh_1$ and $\clh_2$ respectively. Then $V_1\oplus V_2:=((V_1)_{x} \oplus (V_2)_{x})_{x \in P}$ is an isometric representation of $P$ on $\clh_1 \oplus \clh_2$. It is clear that if $V_1$ and $V_2$ are strongly continous then $V_1\oplus V_2$ is strongly continous. Under the identification $\Gamma(\clh_1 \oplus \clh_2) \cong \Gamma(\clh_1) \otimes \Gamma(\clh_2)$, we have the equality $\alpha^{V_1\oplus V_2}=\alpha^{V_1}\otimes \alpha^{V_{2}}$ (See Remark 4.8 of \cite{Anbu}).

What are the examples of isometric representations of $P$ ? One obvious isometric representation is the ``left" regular representation of $P$ on $L^{2}(P)$. One can also consider the ``left" regular representation with multiplicity. A slight generalisation of the above is as follows. Let $A \subset \mathbb{R}^{n}$ be a proper closed subset which is  invariant under translation by elements of $P$ i.e. $A+x \subset A$ for $x \in P$.  Such a subset is called a \emph{$P$-module}. The notion of $P$-modules in the discrete setting was first considered by Salas in \cite{Salas}. Let $k \in \{1,2,\cdots,\} \cup \{\infty\}$ be given and let $\clk$ be a Hilbert space of dimension $k$. Consider the Hilbert space $L^{2}(A,\clk)$. For $x \in P$, let $V_{x}$ be the isometry on $L^{2}(A,\clk)$ defined by the equation: \begin{equation}
\label{isometries}
V_{x}(f)(y):=\begin{cases}
 f(y-x)  & \mbox{ if
} y -x \in A,\cr
   &\cr
    0 &  \mbox{ if } y-x \notin A
         \end{cases}
\end{equation}
for $f \in L^{2}(A,\clk)$. Then $V:=\{V_{x}\}_{x \in P}$ is a strongly continous isometric representation of $P$ on $L^{2}(A,\clk)$. We call $V$ \emph{the isometric representation associated to the $P$-module $A$ of multiplicty $k$}. We call the associated  CCR flow  \emph{the CCR flow corresponding to the $P$-module $A$ of multiplicity $k$} and denote it by $\alpha^{(A,k)}$.  

Let $A$ be a $P$-module and $z \in \mathbb{R}^{d}$ be given. Set $B:=A+z$. Then $B$ is clearly a $P$-module. It is clear that for $k \in \{1,2,\cdots,\}\cup \{\infty\}$, the CCR flow $\alpha^{(A,k)}$ is conjugate to $\alpha^{(B,k)}$. This is due to the fact that the associated isometric represenations are unitarily equivalent. Thus there is no loss of generality in assuming that our $P$-modules contain the origin $0$, which we henceforth assume. The goal of this paper is to prove the following theorem.
\begin{thm}
\label{main theorem}
Let $A_1$ and $A_2$ be $P$-modules and $k_1,k_2 \in \{1,2,\cdots\} \cup \{\infty\}$ be given. Then the following are equivalent.
\begin{enumerate}
\item[(1)] The CCR flow $\alpha^{(A_1,k_1)}$ is  conjugate to $\alpha^{(A_2,k_2)}$.
\item[(2)] The CCR flow $\alpha^{(A_1,k_1)}$ is cocycle conjugate to $\alpha^{(A_{2},k_{2})}$. 
\item[(3)] There exists $z \in \mathbb{R}^{d}$ such that $A_1+z=A_2$ and $k_1=k_2$.
\end{enumerate}
\end{thm}
The implications $(1) \implies (2)$ and $(3) \implies (1)$ are obvious. The difficult part lies in establishing the implication $(2) \implies (3)$. 

The above theorem in the $1$-dimensional case was proved by Arveson. To see this, observe that when $d=1$ i.e. when $\mathbb{R}^{d}$ is $1$-dimensional the only possible choices of $P$ are $[0,\infty)$ or $(-\infty,0]$. Since $[0,\infty)$ and $[-\infty,0)$ are isomorphic, we can assume that $P=[0,\infty)$. Also note that up to a translate, the only $P$-module is $P=[0,\infty)$-itself. Then $\alpha^{([0,\infty),k)}$ is nothing but the usual $1$-dimensional CCR flow of index $k$ and the index is a complete invariant of such CCR flows. For a proof of this well known fact, we refer the reader to \cite{Arveson}. The classification of $1$-dimensional CCR flows relies heavily on the fact that the $1$-dimensional CCR flows have units in abundance. Though we do not need the following fact , we must mention here that the situation in the multi-dimensional case is different and there are not enough units. However another invariant, the gauge group of an $E_0$-semigroup, comes to our rescue. 

Since the classification of the $1$-dimensional CCR flows is complete, we no longer concentrate on the $1$-dimensional case and we assume from now on that the dimension of $\mathbb{R}^{d}$ i.e. $d \geq 2$.  We now explain the techniques behind the proof of Theorem \ref{main theorem}. As already mentioned the gauge group of an $E_0$-semigroup plays a key role in establishing the proof of Theorem \ref{main theorem}. Let us recall the notion of the gauge group associated to an $E_0$-semigroup.

Let $\alpha:=\{\alpha_{x}\}_{x \in P}$ be an $E_0$-semigroup on $B(\clh)$ where $\clh$ is an infinite dimensional separable Hilbert space. An $\alpha$-cocycle is called a \emph{gauge cocycle} if it leaves $\alpha$ invariant. To be precise, let $U:=\{U_{x}\}_{x \in P}$ be an $\alpha$-cocycle. Then $U$ is called a gauge cocycle of $\alpha$ if for $x \in P$ and $A \in B(\clh)$, $U_{x}\alpha_{x}(A)U_{x}^{*}=\alpha_{x}(A)$. The set of gauge cocycles of $\alpha$,  known as  \emph{the gauge group} of $\alpha$ and denoted $G(\alpha)$, is a topological group. For $U:=\{U_{x}\}_{x \in P}$ and $V:=\{V_{x}\}_{x \in P} \in G(\alpha)$ the multiplication $UV$ is given by $UV:=\{U_{x}V_{x}\}_{x \in P}$. The inverse of $U$ is given by $U^{*}:=\{U_{x}^{*}\}_{x \in P}$.  The topology on $G(\alpha)$ is the topology of uniform convergence on compact subsets of $P$, where  the topology that we impose on the unitary group $\mathcal{U}(\clh)$ is the strong operator topology. 

The main result obtained in \cite{Anbu}, which we recall now, is the description of the gauge group of a CCR flow associated to a strongly continous isometric representation which is pure. Let $V:P \to B(\clh)$ be an isometric representation. 
We say that $V$ is \textbf{pure} if for $a \in \Omega$   $V_{ta}^{*}$ converges strongly to zero as $t$ tends to infinity (Recall that the $\Omega$ is the interior of $P$). It is proved in Prop 4.6 of \cite{Anbu} that isometric representations associated to $P$-modules are pure. In what follows, let $V:P \to B(\clh)$ be a strongly continous isometric  representation which is pure. For $a \in P$, we denote the range projection of $V_{a}$ by $E_{a}$. The orthogonal complement of $E_{a}$ i.e. $1-E_{a}$ will be denoted by $E_{a}^{\perp}$. Let $M$ be the commutant of the von Neumann algebra generated by $\{V_{x}:x \in P\}$. Denote the unitary group of $M$ by $\mathcal{U}(M)$. We endow $\mathcal{U}(M)$ with the strong operator topology.

Let $\xi:P \to \clh$ be a map and   denote the image of $x$ under $\xi$ by $\xi_{x}$. We say that $\xi$ is an \emph{additive cocycle of $V$} in case $\xi$ satisfies the following three conditions:  
\begin{enumerate}
\item[(1)] $\xi_{x}+V_{x}\xi_y=\xi_{x+y}$, $x,y \in P$,
\item[(2)]  $V_{x}^{*}\xi_{x}=0$, $x \in P$, and
\item[(3)] $\xi$ is continuous with respect to the norm topology on $\clh$. 
\end{enumerate}
Let $\mathcal{A}(V)$ denote the set of additive cocycles of $V$. We endow $\mathcal{A}(V)$ with the topology of uniform convergence on compact subsets of $P$ where $\clh$ is given the norm topology. The main theorem obtained in \cite{Anbu} (Thm. 7.2) is stated below.
\begin{thm}
\label{Gauge group}
The map \[\mathbb{R}^{d} \times \mathcal{A}(V) \times \mathcal{U}(M) \ni (\lambda,\xi,U) \to \{e^{i\langle \lambda|x\rangle}W(\xi_{x})\Gamma(UE_{x}^{\perp}+E_{x})\}_{x \in P} \in G(\alpha)\] is a homeomorphism.
\end{thm}

Now we explain the contents and the organisation of this paper.

Let $A$ be a $P$-module and $V$ be the isometric representation associated to $A$ of multiplicity $k$ where $k \in \{1,2,\cdots,\}\cup \{\infty\}$. Denote the CCR flow associated to $V$ by $\alpha^{(A,k)}$. In section 2, we show that $V$ admits no non-zero additive cocycle. (Recall that we have assumed that the dimension of the vector space $\mathbb{R}^{d}$ i.e. $d \geq 2$). We achieve this fact by appealing to the theory of distributions. An immediate consequence of the  vanishing of additive cocycles is the fact that the gauge group of $\alpha^{(A,k)}$ is isomorphic to $\mathbb{R}^{d} \times \mathcal{U}(M)$ where $M$ is the commutant of the von Neumann algebra generated by $\{V_{x}:x \in P\}$ and $\mathcal{U}(M)$ is the unitary group of $M$ endowed with the strong operator topology. We must mention here that in \cite{Anbu} this result was obtained for a few examples of $\mathbb{R}_{+}^{2}$-modules using barehand techniques.

In section 3, we compute the commutant $M$ of the von Neumann algebra generated by $\{V_{x}:x \in P\}$.  It is not difficult to see that it suffices to compute the commutant $M$ when the isometric representation $V$ is of multiplicity 1. Let $V:P \to B(L^{2}(A))$ be the isometric representation associated to the $P$-module $A$ of multiplicity $1$ and let $M$ be the commutant of the von Neumann algebra generated by $\{V_{x}:x \in P\}$. Let 
\[
G_{A}:=\{z \in \mathbb{R}^{d}: A+z=A\}.\]
It is clear that $G_{A}$ forms a subgroup of $\mathbb{R}^{d}$. We call $G_{A}$  \emph{the isotropy group} of the $P$-module $A$. For $z \in G_{A}$, let $U_{z}:L^{2}(A) \to L^{2}(A)$ be the unitary defined by the equation 
\[
U_{z}f(x)=f(x-z)\]
for $f \in L^{2}(A)$. 
We show that $M$ is generated by $\{U_{z}:z \in G_{A}\}$. 
Here is where we employ groupoid techniques. The second author has constructed in \cite{Sundar_Ore} a ``universal groupoid" which encodes all isometric representations with commuting range projections. We must mention here that the results obtained in \cite{Sundar_Ore} owes a lot to the papers \cite{Renault_Muhly}, \cite{Nica_WienerHopf} and \cite{Hilgert_Neeb}. In Section 4, we prove our main theorem i.e. Theorem \ref{main theorem}.

Let us end this introduction by thanking a few people who has helped us immensely by  sharing their knowledge on Mathematics whenever we were faced with a difficult problem. First and foremost, we thank R. Srinivasan for introducing us to the beautiful theory of $E_0$-semigroups and also for illuminating discussions on the subject. We would like to thank Murugan  for useful conversations about the symmetric Fock space, the exponential vectors  and the Weyl operators. We thank Prof. Ramadas for directing us towards the theory of distributions in proving Prop.2.4. We thank Prof. V.S. Sunder for sowing the seeds for the proof of Prop. \ref{Sunder}.  
Last but not the least, the second author is hugely indebted to Prof. Renault for passionately sharing his knowledge on groupoids without which this paper would not have materialsed. 

We dedicate this paper in the memory of Prof. Arveson whose ideas have not only inspired us but also many others.

\section{Lack of additive cocycles}
First we collect a few topological and measure theoretical aspects of $P$-modules. Let $A$ be a $P$-module which is  fixed for the rest of this section. Recall that we always assume that $0 \in A$. We denote the interior of $A$ by $Int(A)$ and the boundary of $A$ by $\partial(A)$. For a proof of the following Lemma, we refer the reader to Lemma II.12 of \cite{Hilgert_Neeb}.

\begin{lmma}
We have the following. 
\begin{enumerate}
\item[(1)] The interior of $A$, $Int(A)$ is dense in $A$, and
\item[(2)] The boundary of $A$, $\partial(A)$ has measure zero. 
\end{enumerate}
\end{lmma}
We need the following topological fact in the sequel.
\begin{lmma}
\label{connectedness1}
The Interior of $A$, $Int(A)$ and $A$ are connected. 
\end{lmma}
\textit{Proof.} Note that $A+\Omega \subset Int(A)$.  Let $x,y \in Int(A)$ be given. Since $\Omega$ spans $\mathbb{R}^{d}$, there exists $b,a \in \Omega$ such that $x-y=b-a$ i.e. $x+a=y+b$. Observe that $\{x+ta\}_{t \in [0,1]}$ is a path in $Int(A)$ connecting $x$ and $x+a$. Similarly $\{y+tb\}_{t \in [0,1]}$ is a path  in $Int(A)$ connecting  $y$ and $y+b$. Since $x+a=y+b$, it follows that $x$ is connected to $y$ by a path in $Int(A)$. This proves that $Int(A)$ is path connected and hence connected. Since $\overline{Int(A)}=A$, it follows that $A$ is connected. This completes the proof. \hfill $\Box$

We collect in the following proposition a few facts regarding the topology of $A$ and its boundary $\partial(A)$. For  $E,F \subset \mathbb{R}^{d}$, we denote the complement of $F$ in $E$ by $E\backslash F$. 
\begin{ppsn}
\label{topology of $P$-modules}
Let $a \in \Omega$ be given. We have the following.
\begin{enumerate}
\item[(1)]  The map $\partial(A)\times(0,1) \ni (x,s) \to x+sa \in Int(A)\backslash(A+a)$ is continuous and is onto.
\item[(2)]  The sequence $\{A+na\}_{n \geq 1}$ is a decreasing sequence of closed subsets which decreases to the empty set i.e. $\bigcap_{n=1}^{\infty}(A+na)=\emptyset$. 
\item[(3)] The map $\partial(A) \times (0,\infty) \ni (x,s) \to x+sa \in Int(A)$ is a homemorphism.
\item[(4)] The boundary $\partial(A)$ is connected. 
\item[(5)] The boundary $\partial(A)$ is unbounded.
\item[(6)] The set $Int(A) \backslash(A+a)$ has infinite Lebesgue measure.  
\end{enumerate}
\end{ppsn}
\textit{Proof.} Let $(x,s) \in \partial(A)\times(0,1)$ be given. Since $A+\Omega$ is an open subset of $\mathbb{R}^{d}$ contained in $A$, it follows that $A+\Omega \subset Int(A)$. This implies that $x+sa \in Int(A)$. Now suppose $x+sa \in A+a$. Then there exists $y \in A$ such that $x+sa=y+a$, i.e. $x=y+(1-s)a \in A+\Omega \subset Int(A)$. This contradicts the fact that $x \in \partial(A)$. Hence $x+sa \notin A+a$. Thus we have shown that the map prescribed in $(1)$ is meaningful. The continuity of the prescribed map is obvious. 

Now let $z \in Int(A)\backslash(A+a)$ be given. Let \[E:=\{t \in [0,1]:z-ta \in Int(A)\}.\] Note that $E$ contains $0$ and is an open subset of $[0,1]$.  Denote the supremum of $E$ by $s$. Since $E$ is open in $[0,1]$ and contains $0$, it follows that $s>0$. Note that $z-a \notin A$ and $A$ is a closed subset of $\mathbb{R}^{d}$. Thus for $t$ sufficiently close to $1$, $z-ta \notin A$. This proves that $s<1$. Hence $0<s<1$. As $E$ is open in $[0,1]$, it follows that $s \notin E$ i.e. $z-sa \notin Int(A)$. Choose a sequence $s_n \in E$ such that $s_n \to s$. Then $z-s_na \in Int(A) \subset A$ and $z-s_na\to z-sa$. Since $A$ is closed in $\mathbb{R}^{d}$, it follows that $z-sa \in A$. As a consequence, we have $z-sa \in \partial(A)$. Now note that $z=(z-sa)+sa$. This proves that the map 
\[
\partial(A) \times (0,1) \ni (x,s) \to x+sa \in Int(A)\backslash(A+a)\]
is onto. This proves $(1)$. 

Since $A+P \subset A$, it is clear that $\{A+na\}_{n \geq 1}$ is a decreasing sequence of closed subsets. Suppose $x \in \bigcap_{n=1}^{\infty}(A+na)$. This implies that $x-na \in A$ for every $n \geq 1$. Let $y \in \mathbb{R}^{d}$ be given. By the Archimedean principle there exists a positive integer $n_0$ such that $n_0a-(x-y) \in \Omega$. As $A+\Omega \subset A$, it follows that $y=(x-n_0a)+(n_0a-(x-y))\in A+\Omega \subset A$. This proves that $y \in A$ for every $y \in \mathbb{R}^{d}$ which is a contradiction since $A$ is a proper closed subset of $\mathbb{R}^{d}$. This proves $(2)$. 

The well-definedness of the map in $(3)$ is clear as  $A+\Omega \subset Int(A)$. Let $y \in Int(A)$ be given. By $(2)$, there exists $n \geq 1$ such that $y \notin A+na$. Now by $(1)$, applied to the interior point $na$, it follows that there exists $s \in (0,1)$ and $x \in \partial(A)$ such that $y=x+s(na)$. This proves that the map \[
\partial(A) \times (0,\infty) \ni (x,s) \to x+sa \in Int(A)\]
is onto. Let $x_1,x_2 \in \partial(A)$ and $s_1,s_2 \in (0,\infty)$ be such that $x_1+s_1a=x_2+s_2a$. We claim that $x_1=x_2$ and $s_1=s_2$. It suffices to show that $s_1=s_2$. Suppose not. Without loss of generality, we can assume $s_1>s_2$. Then 
$x_2=x_1+(s_1-s_2)a \in A+\Omega \subset Int(A)$ which contradicts the fact that $x_2 \in \partial(A)$. This proves our claim.  In other words,  the map 
\[
\partial(A) \times (0,\infty) \ni (x,s) \to x+sa \in Int(A)\]
is an injection. It is clear that the above map is continuous. Let $(x_n,s_n) \in \partial(A) \times (0,\infty)$ be a sequence and $(x,s) \in \partial(A) \times (0,\infty)$ be such that $x_n+s_na \to x+sa$. We claim that $x_n \to x$ and $s_n \to s$. It is enough to prove that $s_n \to s$. Suppose $s_n \nrightarrow s$. Then there exists $\epsilon>0$ such that $s_n \notin (s-\epsilon,s+\epsilon)$ for infinitely many $n$. Suppose $s_n \leq s-\epsilon$ for infinitely many $n$. Choose a subsequence $s_{n_k}$ such that $s_{n_k} \in (0,s-\epsilon]$. By passing to a subsequence if necessary we can assume that $s_{n_k}$ converges, say to, $t$. Then $t<s$. Now note that $x_{n_k} \to x+(s-t)a \in A+\Omega \subset Int(A)$.  This is a contradiction since $x_{n_k} \in \partial(A)$ and $\partial(A)$ is a closed subset of $\mathbb{R}^{d}$ which is disjoint from $Int(A)$. 

Now suppose that $s_n \geq s+\epsilon$ for infinitely many $n$. Choose a subsequence $(s_{n_k})$ such that $s_{n_k} \geq s+\epsilon$. Write $s_{n_k}=t_{n_k}+s+\epsilon$ with $t_{n_k} \geq 0$. Note that \begin{align*}
x_{n_k}+t_{n_k}a+\epsilon a &=(x_{n_k}+s_{n_k}a)-sa \\
                                              &\to x
\end{align*} 
But $x_{n_k}+t_{n_{k}}a+\epsilon a \in A+\epsilon a$ and $A+\epsilon a$ is a closed subset of $\mathbb{R}^{d}$. This implies that $x \in A+\epsilon a \subset A+\Omega \subset Int(A)$. This contradicts the fact that $x \in \partial(A)$. These contradictions imply that our assumption $s_n \nrightarrow s$ is wrong and hence $s_n \to s$. Hence the map 
\[
\partial(A) \times (0,\infty) \ni (x,s) \to x+sa \in Int(A)\]
is a homeomorphism. This proves $(3)$. It is immediate that $(3)$ implies $(4)$.

Suppose $\partial(A)$ is bounded. Since $0 \in A$ and $A+\Omega \subset Int(A)$, it follows that $\Omega \subset Int(A)$. Let $b \in \Omega$ be given. By $(3)$, there exists a sequence $(s_n) \in (0,\infty)$ and $x_n \in \partial(A)$ such that $nb-s_na=x_n$. Since $\partial(A)$ is bounded, it follows that $b-\frac{s_n}{n}a=\frac{x_n}{n} \to 0$. In other words, $\frac{s_n}{n}a \to b$. Note that $\{sa: s \geq 0\}$ is a closed subset of $\mathbb{R}^{d}$. Hence there exists $s \geq 0$ such that $b=sa$. Since $\Omega-\Omega=\mathbb{R}^{d}$, it follows that the linear span of $a$ is $\mathbb{R}^{d}$. This implies that $d=1$ which is a contradiction to our assumption that $d \geq 2$. This contradiction implies that $\partial(A)$ is unbounded. 

Let $E:=\{x \in \mathbb{R}^{d}:0< x<a\}$. Then $E$ is a non-empty (as $\frac{a}{2} \in E$) open and bounded set. For a proof of this fact, we refer the reader to the first line of the proof of Prop. I.1.8 in \cite{Faraut}. Let $M:=\sup\{||x||:x \in E\}$ where  the norm on $\mathbb{R}^{d}$ is the usual Euclidean norm. Since $\frac{a}{2} \in E$, it follows that $M \geq \frac{||a||}{2}$. Thus $M>0$. The unboundedness of $\partial(A)$ implies that there exists a sequence $\{x_n\}_{n\geq 1}$ in $\partial(A)$ such that $||x_n-x_m|| \geq 3M$ if $n \neq m$. We claim the following.
\begin{enumerate}
\item[(i)] For $n \geq 1$, $E+x_n \subset Int(A)\backslash(A+a)$, and
\item[(ii)] the sequence $\{E+x_n\}_{n \geq 1}$ forms a disjoint family of non-empty open subsets of $\mathbb{R}^{d}$. 
\end{enumerate}
Let $n\geq 1$ be given. Note that $E \subset \Omega$. Since $\Omega+A \subset Int(A)$, it follows that $E+x_n $ is contained in the interior of $A$. Suppose the intersection $(E+x_n)\cap (A+a)\neq \emptyset$. Then there exists $e \in E$, $y \in A$ such that $e+x_n=y+a$. Hence \[x_n=y+(a-e) \in A+\Omega \subset Int(A).\] This implies that $x_n \in Int(A)$ which contradicts the fact that $x_n \in \partial(A)$. This contradiction implies that the intersection $(E+x_n)\cap (A+a)= \emptyset$.  Hence $E+x_n \subset Int(A)\backslash(A+a)$. This proves $(i)$. 

Let $m,n \geq 1$ be such that $m\neq n$. Suppose that the intersection $(E+x_n)\cap (E+x_m)$ is non-empty. Then there exists $e_1,e_2 \in E$ such that $e_1+x_n=e_2+x_m$. Now  observe that 
\[
3M \leq ||x_n-x_m|| =||e_2-e_1|| \leq ||e_2||+||e_1|| \leq 2M
 \]
which is a contradiction since $M>0$. This proves that the intersection $(E+x_n)\cap (E+x_m)$ is empty. This proves $(ii)$.

Let $\lambda$ be the Lebesgue measure on $\mathbb{R}^{d}$. Since $E$ is a non-empty open subset of $\mathbb{R}^{d}$, it follows that $\lambda(E)>0$. Now calculate as follows to observe that 
\begin{align*}
\infty&=\sum_{n=1}^{\infty}\lambda(E) \\
      &=\sum_{n=1}^{\infty}\lambda(E+x_n)\\
  &=\lambda\big(\coprod_{n=1}^{\infty}E+x_n\big)\\
& \leq \lambda(Int(A)\backslash(A+a)).
\end{align*}
Hence $Int(A)\backslash(A+a)$ has infinite Lebesgue measure. This proves $(6)$ and the proofs are now complete. \hfill $\Box$


The next proposition shows that the isometric representation associated to $A$ of multiplicity $1$ admits no non-zero additive cocycles.
\begin{ppsn}
\label{vanishing of additive cocycles}
Let $V:P \to B(L^{2}(A))$ be the isometric representation associated to the $P$-module $A$ of multiplicity $1$. Suppose that $\{\xi_{x}\}_{x \in P}$ is an additive cocycle of $V$. Then for every $x \in P$, $\xi_{x}=0$. 

\end{ppsn}
\textit{Proof.} Fix $a \in \Omega$. Since $V_{a}^{*}\xi_{a}=0$, it follows that $\xi_{a}(x)=0$ for almost all $x \in A+a$. Without loss of generality, we can assume that $\xi_{a}(x)=0$ for all $x \in A+a$. Also $A\backslash(A+a)$ and $Int(A)\backslash(A+a)$ differ by a set of measure zero. For, the boundary $\partial(A)$ has measure zero. Thus without loss of generality, we can assume that $\xi_{a}(x)=0$ for $x \in \partial(A)$. 

 Let $U:=Int(A)\backslash (A+a)$. By $(1)$ and $(4)$ of Prop.\ref{topology of $P$-modules}, it follows that $U$ is a non-empty open connected subset of $\mathbb{R}^{d}$. Note that the complex conjugate of $\xi_{a}$, i.e. $\overline{\xi_{a}} \in L^{2}(U) \subset L^{1}_{loc}(U)$. Thus, we view $\overline{\xi_{a}}$ as a distribution on $U$. 
 Let $\phi: U \to \mathbb{R}$ be a smooth function such that $supp(\phi)$ is compact and $supp(\phi) \subset U$. Denote the support of $\phi$ by $K$. We view $\phi$ as a smooth function on $\mathbb{R}^{d}$ by declaring its value on the complement of $U$ to be zero. We denote its $i^{th}$ partial derivative of $\phi$ by $\partial_{i}\phi$. For $x \in \mathbb{R}^{d}$, let $\nabla \phi(x)$ be the gradient of $\phi$ i.e. \[\nabla \phi(x)=(\partial_{1}\phi(x),\partial_{2}\phi(x),\cdots,\partial_{n}\phi(x)).\] Let 
 \[
 M:=\sup_{x \in \mathbb{R}^{d}}||\nabla \phi(x)||.
 \]
 Fix $b \in \Omega$  such that $||b||=1$. We claim that there exists $\delta_{0}>0$ such that if $0<t<\delta_{0}$ then  $K \cap (A\backslash(A+tb))=\emptyset$. Suppose not. Then there exists a sequence $(x_n) \in K$ and a sequence of positive real numbers $t_n \to 0$ such that $x_n \in A\backslash(A+t_nb)$. By passing to a subsequence, if necessary, we can assume that $x_n$ converges say to $x \in K$. Note that $K \subset Int(A)$ and $x_n-t_nb \to x \in Int(A)$. Hence eventually $x_n \in Int(A)+t_nb\subset A+t_nb$ which is a contradiction to the fact that $x_{n} \in A \backslash(A+t_nb)$.  This proves our claim. Choose such a $\delta_{0}$. 

Let $\delta>0$  be such that $\delta < \delta_0$ and $K+\overline{B(0,\delta)} \subset U$. Let $L:=K+\overline{B(0,\delta)}$. Note that $L$ is a compact subset of $U$. Let $(t_n)$ be a sequence of positive numbers such that $t_n < \delta$ and $t_n \to 0$. Note that by the  mean value inequality, we have for $x \in \mathbb{R}^{d}$ and $n \geq 1$,  
\begin{equation}
\label{dominated}
\Big|\frac{\phi(x+t_nb)-\phi(x)}{t_{n}}\Big| \leq M 1_{L}(x).
\end{equation}
Note that since $K \cap (A \backslash (A+t_nb))=\emptyset$, the inner product $\langle \xi_{t_nb}|\phi \rangle=0$.  Now calculate as follows to observe that 
\begin{align*}
\int_{U} \frac{\phi(x+t_nb)-\phi(x)}{t_n}\overline{\xi_{a}(x)}dx&= \frac{1}{t_n}\Big(\langle \xi_{a}|V_{t_nb}^{*}\phi \rangle - \langle \xi_{a}|\phi\rangle \Big) \\
                                                                                               & =\frac{1}{t_n}\Big(\langle V_{t_nb}\xi_{a}|\phi \rangle - \langle \xi_{a}|\phi\rangle \Big) \\
                                                                                               &=\frac{1}{t_n}\Big(\langle \xi_{t_nb}+V_{t_nb}\xi_{a}|\phi \rangle - \langle \xi_{a}|\phi\rangle \Big) \\
                                                                                               &=\frac{1}{t_n}\Big(\langle \xi_{t_nb+a}|\phi \rangle - \langle \xi_{a}|\phi\rangle \Big) \\
                                                                                               &=\frac{1}{t_n}\Big(\langle \xi_{a}+V_{a}\xi_{t_nb}|\phi \rangle - \langle \xi_{a}|\phi\rangle \Big) \\
                                                                                               &=\frac{1}{t_n} \langle \xi_{t_nb}|V_{a}^{*}\phi\rangle \\
                                                                                               &=0 ~~(\textrm{~~since $\phi$ vanishes on $A+a$}).
 \end{align*}
 Thus we obtain, for $n \geq 1 $, the equation  \begin{equation}\label{derivative}
 \int_{U} \frac{\phi(x+t_nb)-\phi(x)}{t_n}\overline{\xi_{a}(x)}dx=0 
 \end{equation}
 For $x \in U$ and $n \geq 1$, Eq.\ref{dominated} implies that \[\Big|\Big(\frac{\phi(x+t_nb)-\phi(x)}{t_n}\Big)\overline{\xi_{a}(x)}\Big| \leq M1_{L}(x)|\xi_{a}(x)|.\] The function $U \ni x \to 1_{L}(x)|\xi_{a}(x)| \in [0,\infty)$ is integrable.  Thus letting $n \to \infty$ in Eq.\ref{derivative} and applying the dominated convergence theorem, we obtain 
 \[
 \int_{U} \langle \nabla \phi(x)|b \rangle \overline{\xi_{a}(x)}=0.
 \]
Since $t\Omega=\Omega$ for every $t>0$, it follows that for every $b \in \Omega$, 
\[
 \int_{U} \langle \nabla \phi(x)|b \rangle \overline{\xi_{a}(x)}=0.
 \]
 Since $\Omega$ is spanning, it follows that for every $z \in \mathbb{R}^{d}$, $\int_{U}\langle \nabla \phi(x)|z \rangle \overline{\xi_{a}(x)}=0.$ Let $e_1,e_2,\cdots,e_d$ be the standard orthonormal basis of $\mathbb{R}^{d}$. Then for every $i=1,2,\cdots,d$,
 \[
 \int_{U}\partial_{i}\phi(x) \overline{\xi_{a}(x)} dx= \int_{U}\langle \nabla \phi(x)|e_i \rangle \overline{\xi_{a}(x)}=0.
 \] 
Hence each partial derivative of $\overline{\xi_{a}}$, in the distribution sense, vanishes. By Theorem 6.3-4 of \cite{Philippe}, it follows that there exists a complex number $c_{a}$ such that $\xi_{a}(x)=c_{a}$ for almost all $x \in U$. Since $\xi_{a} \in L^{2}(A)$ and $U$ has infinite measure by Prop.\ref{topology of $P$-modules}, it follows that $c_a=0$. Hence $\xi_{a}=0$ for every $a \in \Omega$. The density of $\Omega$ in $P$ and the continuity of the map $P \ni \xi \to \xi_{x} \in L^{2}(A)$ implies that $\xi_{x}=0$ for every $x \in P$. This completes the proof. \hfill $\Box$ 


\begin{rmrk}
\label{additive cocycles}
Let $k \in \{1,2,\cdots\}\cup \{\infty\}$ be given. For each $i$, let $\clh_i$ be a Hilbert space. Denote the direct sum $\bigoplus \clh_i$ by $\clh$. For a vector $\xi \in \clh$, we denote its $i^{th}$-component by $\xi^{i}$. For each $i$, let $V^{i}:=\{V^{i}_{x}\}_{x \in P}$ be  an isometric representation of $P$ on $\clh_{i}$. Let $V:=\bigoplus V^{i}$ be the direct sum. Then cleary $V$ is an isometric representation of $P$ on $\clh$. If each $V^{i}$ is strongly continous then $V$ is strongly continous.

If $\xi:=\{\xi_{x}\}_{x \in P}$ is an additive cocycle of $V$, then $\xi^{i}:=\{\xi_{x}^{i}\}_{x \in P}$ is an additive cocycle of $V^{i}$ for each $i$. Thus if each $V^{i}$ admits no non-trivial additive cocycle then $V$ admits no non-trivial additive cocycle. 
\end{rmrk}

Now the following is an immediate corollary of Prop.\ref{vanishing of additive cocycles} and Remark \ref{additive cocycles}.

\begin{crlre}
\label{vanishing of additive cocycles 1}
Let $k \in \{1,2,\cdots,\}\cup\{\infty\}$ and $\mathcal{K}$ be a Hilbert space of dimension $k$. Let $V:P \to B(L^{2}(A)\otimes \clk)$ be the isometric representation associated to the $P$-module $A$ of multiplicity $k$. Suppose that $\xi:=\{\xi_{x}\}_{x \in P}$ is an additive cocycle of $V$. Then $\xi_{x}=0$ for every $x \in P$. 
\end{crlre}

Corollary \ref{vanishing of additive cocycles 1} and Theorem \ref{Gauge group} leads to the next Theorem. 

 
Let $k \in \{1,2,\cdots,\}\cup\{\infty\}$ and $\clk$ be a Hilbert space of dimension $k$. Let $V$ be the isometric representation of $P$ on the Hilbert space $L^{2}(A)\otimes \clk$  associated to the $P$-module $A$ of multiplicity $k$. For $x \in P$, we denote the range projection of $V_{x}$ by $E_{x}$. The orthogonal complement of $E_{x}$ i.e. $1-E_{x}$ will be denoted by $E_{x}^{\perp}$.
Denote the commutant of the von Neumann algebra generated by $\{V_{x}:x \in P\}$ by $M$. Denote the unitary group of $M$ by $\mathcal{U}(M)$. We endow $\mathcal{U}(M)$ with the strong operator topology.  Let $\alpha^{(A,k)}$ be the CCR flow associated to the isometric representation $V$ and denote the gauge group of $\alpha^{(A,k)}$ by $G(\alpha^{(A,k)})$.

\begin{thm}
\label{description of the gauge group}
With the foregoing notation, the map 
\[\mathbb{R}^{d}  \times \mathcal{U}(M) \ni (\lambda,U) \to \{e^{i\langle \lambda|x\rangle}\Gamma(UE_{x}^{\perp}+E_{x})\}_{x \in P} \in G(\alpha^{(A,k)})\]
  is a homeomorphism where  the topology on $\mathbb{R}^{d} \times \mathcal{U}(M)$ is the product topology.
\end{thm}

\section{The commutant calculation}
Let $A$ be a $P$-module, $k \in \{1,2,\cdots\}\cup\{\infty\}$ and $\clk$ be a Hilbert space of dimension $k$.  Let $V:P \to B(L^{2}(A) \otimes \clk)$ be the isometric representation associated to $A$ of multiplicity $k$. The goal of this section is to compute the commutant of the von Neumann algebra generated by $\{V_{x}:x \in P\}$. First we compute the commutant when the multiplicity $k=1$. This relies heavily on the groupoid approach developed in \cite{Sundar_Ore} to study $C^{*}$-algebras arising out of Ore semigroup actions. 

We must mention here that the results obtained in \cite{Sundar_Ore} are due to the deep insight of Muhly and Renault in using groupoids to understand the Wiener-Hopf $C^{*}$-algebras. This is achieved in their seminal paper \cite{Renault_Muhly}. This view was further developed by Nica in \cite{Nica_WienerHopf} and Hilgert and Neeb in \cite{Hilgert_Neeb}. The results obtained in \cite{Sundar_Ore} also owes a lot to \cite{Hilgert_Neeb}. For completeness, we review the basics of groupoid $C^{*}$-algebras.  For a quick introduction to the theory of groupoids and the associated $C^{*}$-algebras, we either recommend the first two sections of \cite{Khoshkam_Skandalis} or the second section of \cite{Renault_Muhly}. For a more detailed study of groupoids, we refer the reader to  \cite{Renault1}. We recall  here the basics of $C^{*}$-algebras associated to a topological groupoid.

Let $\clg$ be a topological groupoid with a left Haar system. We assume that $\clg$ is locally compact, Hausdorff and second countable. The unit space of $\clg$ will be denoted by $\clg^{(0)}$ and let $r,s:\clg \to \clg^{(0)}$ be the range and source maps. For $x \in \clg^{(0)}$, let $\mathcal{G}^{(x)}=r^{-1}(x)$. Fix a left Haar system $(\lambda^{(x)})_{x \in \clg^{(0)}}$. Let $C_{c}(\clg)$ be the space of continuous complex valued functions defined on $\clg$. The space $C_{c}(\clg)$ forms a $*$-algebra where the multiplication and the involution are defined by the following formulas
\begin{align*}
f*g(\gamma)&=\int f(\eta)g(\eta^{-1}\gamma)d\lambda^{(r(\gamma))}(\eta)\\
f^{*}(\gamma)&=\overline{f(\gamma^{-1})}
\end{align*}
for $f,g \in C_{c}(\clg)$. We obtain bounded representations of the $*$-algebra $C_{c}(\clg)$ as follows. Fix a point $x \in \clg^{(0)}$. Consider the Hilbert space $L^{2}(\clg^{(x)},\lambda^{(x)})$. For $f \in C_{c}(\clg)$ and $\xi \in L^{2}(\clg^{(x)},\lambda^{(x)})$, let $\pi_{x}(f)\xi \in L^{2}(\clg^{(x)},\lambda^{(x)})$ be defined by the formula
\[
(\pi_{x}(f)\xi)(\gamma):=\int f(\gamma^{-1}\gamma_1)\xi(\gamma_1)d\lambda^{(x)}(\gamma_1)\]
for $\gamma \in \clg^{(x)}$.
Then $\pi_{x}:C_{c}(\clg)\to B(L^{2}(\clg^{(x)},\lambda^{(x)})$ is a non-degenerate $*$-representation. Moreover $\pi_{x}$ is continuous when $C_{c}(\clg)$ is given the inductive limit topology. For $f \in C_{c}(\clg)$, let 
\[
||f||_{red}:=\sup_{x \in \clg^{(0)}}||\pi_{x}(f)||.\]
Then $||~||_{red}$ is well defined and is a $C^{*}$-norm on $C_{c}(\clg)$. The completion of $C_{c}(\clg)$ with respect to the norm $||~||_{red}$ is called \emph{the reduced $C^{*}$-algebra of $\clg$} and denoted $C_{red}^{*}(\clg)$. There is also a universal $C^{*}$-algebra associated to $\clg$ and denoted $C^{*}(\clg)$. However we do not need the universal one as the groupoids that we consider are amenable and for amenable groupoids $C_{red}^{*}(\clg)$ and $C^{*}(\clg)$ coincide. Fix $x \in \clg^{(0)}$. We denote the extension of $\pi_{x}$ to $C^{*}(\clg)$ by $\pi_{x}$ itself. The representation $\pi_{x}$ is called the representation of $C^{*}(\clg)$ \emph{induced at the point} $x$. 

Let $\gamma \in \clg$ be such that $s(\gamma)=x$ and $r(\gamma)=y$. Let  $U_{\gamma}:L^{2}(\mathcal{G}^{(y)},\lambda^{(y)}) \to L^{2}(\clg^{(x)},\lambda^{(x)})$  be defined by the formula
\[
U_{\gamma}\xi(\gamma_1)=\xi(\gamma \gamma_1)
\]
for $\xi \in L^{2}(\clg^{(y)},\lambda^{(y)})$. The fact that $(\lambda^{(x)})_{x \in \clg^{(0)}}$ is a left Haar system implies that $U_{\gamma}$ is a unitary. Moreover it is routine to verify that $U_{\gamma}$ intertwines the representations $\pi_{x}$ and $\pi_{y}$, i.e. for $f \in C_{c}(\clg)$, $U_{\gamma}\pi_{y}(f)=\pi_{x}(f)U_{\gamma}$.

We need the following two facts.
\begin{enumerate}
\item[(1)] Let $x \in \clg^{(0)}$ be given. Denote the isotropy group at $x$ by $\clg^{x}_{x}$ i.e. \[\clg^{x}_{x}:=\{\gamma \in \clg: r(\gamma)=s(\gamma)=x\}.\] Note that $\clg^{x}_{x}$ is a group. The commutant of $\{\pi_{x}(f):f \in C_{c}(\clg)\}$ is generated by $\{U_{\gamma}:\gamma \in \clg^{x}_{x}\}$. 
\item[(2)] For $x,y \in \clg^{(0)}$, the representations $\pi_{x}$ and $\pi_{y}$ are non-disjoint if and only if there exists $\gamma \in \clg$ such that $s(\gamma)=x$ and $r(\gamma)=y$. 
\end{enumerate}
Connes proved the above two facts in his paper \cite{Connes_Sur_la}.  In the appendix,  we offer a proof for $(1)$ and $(2)$   for Deaconu-Renault groupoids considered by the second author  and Renault in \cite{Jean_Sundar} which is all we need.  We  believe that the appendix is interesting on its own right as the proof uses notions like groupoid equivalence and Rieffel's notion of strong Morita equivalence.

Let us recall the Deaconu-Renault groupoid considered  in \cite{Jean_Sundar}.  Let $X$ be a compact metric space. By \emph{an action} of $P$ on $X$, we mean a continuous map $X \times P \ni (x,t) \to x+t \in X$ such that $x+0=x$ and $(x+s)+t=x+(s+t)$ for $x \in X$ and $s,t \in P$. We assume that the action of $P$ on $X$ is injective i.e. for every $t \in P$, the map $X \ni x \to x+t \in X$ is injective. Let 
\[
X \rtimes P:=\{(x,t,y) \in X \times \mathbb{R}^{d} \times X: \exists ~r,s \in P~\textrm{such that}~t=r-s \textrm{~and~}x+r=y+s\}.
\]
The set $X \rtimes P$ has a groupoid structure where the groupoid multiplication and the inversion are given by 
\begin{align*}
(x,s,y)(y,t,z)&=(x,s+t,z), ~~\textrm{and} \\
(x,s,y)^{-1}&=(y,-s,x).
\end{align*}
We call $X \rtimes P$ the \emph{Deaconu-Renault groupoid} determined by the action of $P$ on $X$.
The set $X \rtimes P$ is a closed subset of $X \times \mathbb{R}^{d} \times X$. When endowed with the subspace topology, $X \rtimes P$ becomes a topological groupoid. The map $
X \rtimes P \ni (x,t,y) \to (x,t) \in X \times \mathbb{R}^{d}$
is an embedding and the range of the prescribed map is a closed subset of $X \times \mathbb{R}^{d}$. From here on, we always consider $X \rtimes P$ as a subspace of $X \times \mathbb{R}^{d}$. 

For $x \in X$, let \[
Q_{x}:=\{t \in \mathbb{R}^{d}: (x,t) \in X \rtimes P\}.\] Note that for $x \in X$, $Q_{x}$ is a closed subset of $\mathbb{R}^{d}$ containing the origin $0$ and $Q_{x}+P \subset Q_{x}$. By Lemma 4.1 of \cite{Jean_Sundar}, it follows that $Int(Q_{x})$ is dense in $Q_{x}$ and the boundary $\partial(Q_{x})$ has Lebesgue measure zero. 

For $x \in X$, let $\lambda^{(x)}$ be the measure on $X \rtimes P$ defined by the following formula: For $f \in C_{c}(X \rtimes P)$, 
\begin{equation}
\label{Haar system}
\int f d\lambda^{(x)}=\int f(x,t)1_{Q_{x}}(t)dt.
\end{equation}
Here $dt$ denotes the usual Lebesgue measure on $\mathbb{R}^{d}$. The groupoid $X \rtimes P$ admits a Haar system if and only if the map $X \times \Omega \ni (x,s) \to x+s \in X$ is open. In such a  case, $(\lambda^{(x)})_{x \in X}$ forms a left Haar system. When $X \rtimes P$ admits a Haar system, we use only the Haar system described above.

Assume that $X \rtimes P$ has a Haar system. Then the action of $P$ on $X$ can be dilated to an action of $\mathbb{R}^{d}$ on a locally compact space $Y$. More precisely, there exists a locally compact Hausdorff space $Y$, an action of $\mathbb{R}^{d}$ on $Y$, $Y \times \mathbb{R}^{d} \ni (y,t) \to y+t$, an embedding $i:X \to Y$ such that 
\begin{enumerate}
\item[(1)] the embedding $i:X \to Y$ is $P$-equivariant,
\item[(2)] the set $X_{0}:=i(X)+\Omega$ is open in $Y$, and
\item[(3)] the set $Y=\displaystyle \bigcup_{t \in P}(i(X)-t)=\bigcup_{t \in \Omega}(X_{0}-t)$.
\end{enumerate}
The space $Y$ is uniquely determined by conditions $(1)$, $(2)$ and $(3)$ up to an $\mathbb{R}^{d}$-equivariant isomorphism. We suppress the notation $i$ and simply identify $X$ as a subspace of $Y$. We call the pair $(Y,\mathbb{R}^{d})$ the dilation associated to the pair $(X,P)$. Moreover the groupoid $X \rtimes P$ is merely the reduction of the transformation groupoid $Y \rtimes \mathbb{R}^{d}$ onto $X$, i.e. $X \rtimes P:=(Y \rtimes \mathbb{R}^{d})|_{X}$. Also the groupoids $X \rtimes P$ and $Y \rtimes \mathbb{R}^{d}$ are equivalent in the sense of \cite{MRW}. Since $Y \rtimes \mathbb{R}^{d}$ is amenable, it follows from Theorem 2.2.17 of \cite{Claire_Renault} that $X \rtimes P$ is amenable. 
For proofs and details of the facts about Deaconu-Renault groupoids (mentioned in the previous paragraphs), we refer the reader to \cite{Jean_Sundar}.

Let $X$ be a compact metric space and $X \times P \ni (x,t) \to x+t \in X$ be an action of $P$ on $X$. Assume that $X\rtimes P$ has a Haar system. Denote the dilation associated to $(X,P)$ by $(Y,\mathbb{R}^{d})$. Let $\clg:=X \rtimes P$ and $\clh:=Y \rtimes \mathbb{R}^{d}$. The range and source maps of both $\clg$ and $\clh$ will be denoted by $r$ and $s$ respectively. Fix $x \in X$ and let \[Q_{x}:=\{t \in \mathbb{R}^{d}:x+t \in X\}.\] Note that $\clg^{(x)}:=r^{-1}(x):=\{x\} \times Q_{x}$. Thus the Hilbert space $L^{2}(\clg^{(x)},\lambda^{(x)})$ can be identified with $L^{2}(Q_{x})=L^{2}(Q_{x},dt)$ where $dt$ denotes the Lebesgue measure on $\mathbb{R}^{d}$. Let $\clg^{x}_{x}$  be the isotropy group of $\clg$ at $x$ i.e. $\clg^{x}_{x}:=\{t \in \mathbb{R}^{d}:x+t=x\}$. Note that for $t \in \mathbb{R}^{d}$, $t \in \clg^{x}_{x}$ if and only if $(x,t,x) \in X \rtimes P$. For $s \in \clg^{x}_{x}$, let $U_{s}$ be the unitary on $L^{2}(Q_{x})$ defined by the following formula
\[
(U_{s}\xi)(t)=\xi(t-s)\]
for $\xi \in L^{2}(Q_{x})$. 

\begin{thm}
\label{main theorem of the appendix}
With the foregoing notation, we have the following.
\begin{enumerate}
\item[(1)] For $x \in X$, let $\pi_{x}$ be the representation of $C^{*}(\clg)$ induced at $x$. Then the commutant of $\{\pi_{x}(f):f \in C_{c}(\clg)\}$ is the von Neumann algebra generated by $\{U_{s}:s \in \clg^{x}_{x}\}$. 
\item[(2)] Let $x,y \in X$ be given. Then $\pi_{x}$ and $\pi_{y}$ are non-disjoint if and only if there exists $t \in \mathbb{R}^{d}$ such that $x+t=y$ i.e. $(x,t,y) \in \clg$. 
\end{enumerate}
\end{thm}
 We provide a proof of the above theorem in the appendix. 
 
 Let $V:P \to B(\clh)$ be a strongly continous isometric representation with commuting range projections. More precisely, let $E_{x}:=V_{x}V_{x}^{*}$ for $x \in P$. We say that $V$ has \emph{commuting range projections} if $\{E_{x}:x \in P\}$ is a commuting family of projections.  For $z \in \mathbb{R}^{d}$, write $z=x-y$ with $x,y \in P$ and let $W_{z}:=V_{y}^{*}V_{x}$. Then $W_{z}$ is well-defined and is a partial isometry for every $z \in \mathbb{R}^{d}$. Also $\{W_{z}\}_{z \in \mathbb{R}^{d}}$ forms a strongly continuous family of partial isometries. We refer the reader to Prop.3.4 of \cite{Sundar_Ore} for proofs of the above mentioned facts. For $f \in C_{c}(\mathbb{R}^{d})$, let \[W_{f}:= \int f(z)W_{z}dz.\] For $f \in C_{c}(\mathbb{R}^{d})$, $W_{f}$ is called  \emph{the Wiener-Hopf operator} with symbol $f$. 
 
 \begin{lmma}
 \label{von Neumann algebra generated by Wiener-Hopf operators}
 With the foregoing notation, the von Neumann algebra generated by  the set of Wiener-Hopf operators $\{W_{f}:f \in C_{c}(\mathbb{R}^{d})\}$ and the von Neumann algebra generated by $\{V_{x}:x \in P\}$ coincide. 
  \end{lmma}
  \textit{Proof.} It is clear that the von Neumann algebra generated by $\{W_{f}:f \in C_{c}(\mathbb{R}^{d})\}$ is contained in the von Neumann algebra generated by $\{V_{x}:x \in P\}$. Let $z_0 \in \mathbb{R}^{d}$ be given. For $n \geq 1$, let $B(z_0,\frac{1}{n})$ be the open ball centred at $z_0$ and of radius $\frac{1}{n}$. For $n \geq 1$, choose a function $f_n \in C_{c}(\mathbb{R}^{d})$ such that $f_n \geq 0$,  $\int f_{n}(z)dz=1$ and $supp(f_n) \subset B(z_0,\frac{1}{n})$. We claim that $\int f_{n}(z)W_{z}dz \to W_{z_0}$ weakly. 
  
  Let $\xi,\eta \in \clh$ and $\epsilon>0$ be given. Since $\{W_{z}\}_{z \in \mathbb{R}^{d}}$ is strongly continuous, it follows that there exists  $N \geq 1$ such that $|\langle W_{z}\xi|\eta \rangle-\langle W_{z_0}\xi|\eta\rangle| \leq \epsilon$ for every $z \in B(z_0,\frac{1}{N})$. Let $n \geq N$ be given. Calculate as follows to observe that 
  \begin{align*}
  \big| \langle \big(\int f_{n}(z)W_{z}dz\big)\xi|\eta \rangle-\langle W_{z_0}\xi|\eta \rangle \big|&=\big|\int  f_{n}(z)\langle W_{z}\xi|\eta\rangle dz-\int f_{n}(z)\langle W_{z_0}\xi|\eta\rangle dz \big|\\
                                                                                                                                     &=\big| \int f_{n}(z)(\langle W_{z}\xi|\eta\rangle-\langle W_{z_0}\xi|\eta \rangle)dz \big| \\
                                                                                                                                     &\leq \displaystyle \int_{z \in B(z_0,\frac{1}{N})} f_{n}(z)\big|\langle W_{z}\xi|\eta \rangle -\langle W_{z_0}\xi|\eta \rangle\big|dz\\
                                                                                                                                     & \leq \epsilon \int f_{n}(z)dz \\
                                                                                                                                     &\leq \epsilon
   \end{align*}
   This proves that $\int f_{n}(z)W_{z}dz \to W_{z_0}$. Now it is immediate that the von Neumann algebra generated by $\{V_{x}:x \in P\}$ is contained in the von Neumann algebra generated by $\{W_{f}:f \in C_{c}(\mathbb{R}^{d})\}$. This completes the proof. \hfill $\Box$
   
   Next we recall the universal groupoid constructed in \cite{Sundar_Ore}. Denote the set of closed subsets of $\mathbb{R}^{d}$ by $\mathcal{C}(\mathbb{R}^{d})$.  Let \[X_{u}:=\{A \in \mathcal{C}(\mathbb{R}^{d}):0 \in A, -P+A \subset A\}.\] 
     Consider $L^{\infty}(\mathbb{R}^{d})$ as the dual of $L^{1}(\mathbb{R}^{d})$ and endow $L^{\infty}(\mathbb{R}^{d})$ with the weak$^{*}$-topology. The map $X_{u} \ni A \to 1_{A} \in L^{\infty}(\mathbb{R}^{d})$ is injective. Via this injection, we view $X_{u}$ as a subset of $L^{\infty}(\mathbb{R}^{d})$ and endow $X_{u}$ with the subspace topology inherited from the weak$^{*}$-topology on $L^{\infty}(\mathbb{R}^{d})$. The space $X_{u}$ is a compact metric space. The map \[X_{u} \times P \ni (A,x) \to A+x \in X_{u}\]provides us with an injective action of $P$ on $X_{u}$. Moreover the map $X_{u} \times \Omega \ni (A,x) \to A+x \in X_{u}$ is open. See Prop.4.4 and Remark 4.5 of \cite{Sundar_Ore} for a proof of this fact. Consequently, the Deaconu-Renault groupoid $X_{u} \rtimes P$ has a Haar system. We denote the Deaconu-Renault groupoid $X_{u} \rtimes P$ by $\clg_{u}$. The range and source maps of $\clg_u$ will be denoted by $r$ and $s$ respectively. We need the following two facts about the groupoid $\clg_{u}$. For proofs, we refer the reader to Remark 4.5 and Prop.2.1 of \cite{Sundar_Ore}.
     \begin{enumerate}
     \item[(1)] For $A \in X_{u}$, let $Q_{A}:=\{z \in \mathbb{R}^{d}: (A,z) \in X_{u} \rtimes P\}$. Then $Q_{A}=-A$ for every $A \in X_{u}$. 
               \item[(2)] For $f \in C_{c}(\mathbb{R}^{d})$, let $\widetilde{f} \in C_{c}(\clg_u)$  be defined by the equation
     \begin{equation}
     \label{f}
     \widetilde{f}(A,z):=f(z)
     \end{equation}
     for $(A,z) \in \clg_{u}$. Then $C^{*}(\clg_u)$ is generated by $\{\widetilde{f}:f \in C_{c}(\mathbb{R}^{d})\}$. 
     \end{enumerate}
     Let $(\lambda^{(A)})_{A \in X_{u}}$ be the Haar system on $\clg_u$ defined by the equation \ref{Haar system}. Fix $A \in X_{u}$. Then $\clg_{u}^{A}:=r^{-1}(A)=\{A\} \times Q_{A}=\{A\} \times -A$. We identify $L^{2}(\clg_{u}^{A},\lambda^{A})$ with $L^{2}(-A)$. 
     
 Let $A$ be a $P$-module and $V:P \to B(L^{2}(A))$ be the isometric representation associated to $A$ of multiplicity $1$. Let $\{W_{z}\}_{z \in \mathbb{R}^{d}}$ be the partial isometries, described in the paragraph following Theorem \ref{main theorem of the appendix}, associated to the isometric representation $V$. Note that $-A \in X_{u}$. Denote the representation of $C^{*}(\clg_{u})$ induced at $-A$ by $\pi_{A}$. 
 
 \begin{ppsn}
 \label{Induced representation}.
 With the foregoing notation, we have 
 \[
 \pi_{A}(\widetilde{f})=\int f(-z)W_{z}dz
 \]
 for every $f \in C_{c}(\mathbb{R}^{d})$. Here for $f \in C_{c}(\mathbb{R}^{d})$, $\widetilde{f} \in C_{c}(\clg_u)$ is as defined in Equation \ref{f}. 
 
 \end{ppsn}
 We omit the proof of the above proposition as it is  similar to the calculations carried out in the two paragraphs following Remark 5.3 of \cite{Jean_Sundar}. 
 
 Fix a $P$-module say $A$ for the rest of this section. Let $V:P \to B(L^{2}(A))$ be the isometric representation associated to $A$ of multiplicity $1$. Let $G_{A}$ be the isotropy group of $A$, i.e. \[G_{A}:=\{z \in \mathbb{R}^{d}:A+z=A\}.\] For $z \in G_{A}$, let $U_{z}$ be the unitary defined on $L^{2}(A)$ by the equation
\begin{equation}
\label{unitaries}
U_{z}f(x):=f(x-z)
\end{equation}
for $f \in L^{2}(A)$. 

The following corollary is an immediate consequence of Theorem \ref{main theorem of the appendix},  Prop. \ref{Induced representation}, Lemma \ref{von Neumann algebra generated by Wiener-Hopf operators} and the fact that $\{\widetilde{f} \in C_{c}(\clg_u):f \in C_{c}(\mathbb{R}^{d})\}$ generates $C^{*}(\clg_u)$. 

\begin{crlre}
\label{commutant 1}
With the foregoing notation, we have that the commutant of the von Neumann algebra generated by $\{V_{x}:x \in P\}$ coincides with the von Neumann algebra generated by $\{U_{z}:z \in G_{A}\}$. 
\end{crlre}

Let $k \in \{1,2,\cdots,\}\cup\{\infty\}$ be given and let $\clk$ be a Hilbert space of dimension $k$. Let $V:P \to B(L^{2}(A) \otimes \clk)$ be the isometric representation associated to $A$ of multiplicity $k$. Denote the isometric representation associated to $A$ of multiplicity $1$ by $\widetilde{V}$. Then it is clear that for $x \in P$, $V_{x}=\widetilde{V}_{x} \otimes 1$. Let $N$ be the von Neumann algebra on $L^{2}(A)$ generated by $\{\widetilde{V}_{x}:x \in P\}$. Denote the commutant of $N$ by $M$ i.e. $M=N^{'}$. Corollary \ref{commutant 1} implies that $M$ is generated by $\{U_{z}:z \in G_{A}\}$. In particular, $M$ is abelian. Denote the CCR flow associated to the isometric representation $V$ by $\alpha^{(A,k)}$. Let $G(\alpha^{(A,k)})$ denote the gauge group of $\alpha^{(A,k)}$. The following corollary is now immediate.

\begin{crlre}
\label{commutant 2}
With the foregoing notation, the commutant of the von Neumann algebra generated by $\{V_{x}:x \in P\}$ is $M \otimes B(\clk)$. The gauge group of $\alpha^{(A,k)}$ i.e. $G(\alpha^{(A,k)})$ is isomorphic to $\mathbb{R}^{d} \times \mathcal{U}(M \otimes B(\clk))$, where $\mathcal{U}(M \otimes B(\clk))$ is the unitary group of $M \otimes B(\clk)$ endowed with the strong operator topology. 
\end{crlre}
\textit{Proof.} The von Neumann algebra generated by $\{V_{x}:x \in P\}$ is $N \otimes 1:=\{T \otimes 1: T \in N\}$. Hence the commutant of $N \otimes 1$ is $M \otimes B(\clk)$. The fact that $G(\alpha^{(A,k)})$ is isomorphic to $\mathbb{R}^{d} \times \mathcal{U}(M \otimes B(\clk))$ follows from Theorem \ref{description of the gauge group}. This completes the proof. \hfill $\Box$
 
 \section{Proof of the main theorem}
 In this section, we prove  \textbf{Theorem \ref{main theorem}}. We need a basic fact regarding the representation theory of the unitary group of an $n$-dimensional Hilbert space. We start  with a  combinatorial lemma. 
 Fix $\ell \geq 2$. Denote the permutation group on $\{1,2,\cdots,\ell\}$ by $S_{\ell}$. For $i,j \in \{1,2,\cdots,\ell\}$ and $i\neq j$, the permutation that interchanges $i$ and $j$ and leaves the rest fixed will be denoted by $(i,j)$. For $\sigma \in S_{\ell-1}$, let $\widehat{\sigma} \in S_{\ell}$ be defined by $\widehat{\sigma}(i)=\sigma(i)$ for $1\leq i \leq \ell-1$ and $\widehat{\sigma}(\ell)=\ell$. Via the embedding $S_{\ell-1} \ni \sigma \to \widehat{\sigma} \in S_{\ell}$, we view $S_{\ell-1}$ as a subgroup of $S_{\ell}$. For $m:=(m_1,m_2,\cdots,m_\ell) \in \mathbb{Z}^{\ell}$ and $\sigma \in S_{\ell}$, let \[m_{\sigma}:=(m_{\sigma(1)},m_{\sigma(2)},\cdots,m_{\sigma(\ell)}).\] 
 
 \begin{lmma}
 \label{combinatorial lemma}
 Let $\ell \geq 2$ and $m:=(m_1,m_2,\cdots,m_\ell) \in \mathbb{Z}^{\ell}$ be given. Suppose that there exists $i,j \in \{1,2,\cdots,\ell\}$ such that $i\neq j$ and $m_i \neq m_j$. Then the cardinality of the set $\{m_{\sigma}:\sigma \in S_{\ell}\}$ is at least $\ell$. 
  \end{lmma}
  \textit{Proof.} We prove this by induction on $\ell$. The base case i.e. when $\ell=2$ is clearly true. Fix $\ell \geq 3$ and assume that the conclusion of the Lemma holds for $\ell-1$. Choose $i,j \in \{1,2,\cdots,\ell\}$ such that $i\neq j$ and $m_i \neq m_j$. Replacing $m$ by $m_{\sigma}$ for a suitable $\sigma$ if necessary, we can without loss of generality assume that $i=1$. 
  
  \textit{Case 1:} $1 \leq j \leq \ell-1$. Then by the induction hypothesis, the cardinality of the set $\{m_{\widehat{\sigma}}: \sigma \in S_{\ell-1}\}$ is at least $\ell-1$. Since $m_1 \neq m_j$ either $m_1 \neq m_{\ell}$ or $m_{j} \neq m_{\ell}$. Hence there exists $i \in \{1,j\}$ such that $m_i \neq m_\ell$. Then $\{m_{(i,\ell)}\}$ is disjoint from $\{m_{\widehat{\sigma}}:\sigma \in S_{\ell-1}\}$. This proves that the cardinality of the set $\{m_{\sigma}: \sigma \in S_{\ell}\}$ is at least $\ell$. 
  
  \textit{Case 2:} $j=\ell$. If there exists $k \in \{1,2,\cdots,\ell-1\}$ such that $m_1 \neq m_k$ then by Case 1, we have that the cardinality of the set $\{m_{\sigma}: \sigma \in S_{\ell}\}$ is at least $\ell$. Now assume that $m_1=m_k$ for every $k \in \{1,2,\cdots,\ell-1\}$. Then note that the cardinality of the set $\{m_{(i,\ell)}: 1 \leq i \leq \ell\}$ is $\ell$ and hence the cardinality of the set $\{m_{\sigma}: \sigma \in S_{\ell}\}$ is at least $\ell$. This completes the proof. \hfill $\Box$
 
 Consider the $n$-dimensional Hilbert space $\mathbb{C}^{n}$ with the usual Euclidean inner product. For $i=1,2,\cdots,n$, let $e_{i} \in \mathbb{C}^{n}$ be the vector which has $1$ in the $i^{th}$-coordinate and zero elsewhere. Denote the unitary group of $\mathbb{C}^{n}$ endowed with the norm topology by $U(n)$. The special unitary group i.e. the set of unitary operators with determinant one will be denoted by $SU(n)$. Let $\mathbb{T}^{n}$ be the subgroup of $U(n)$ consisting of diagonal matrices. For $\lambda \in \mathbb{T}^{n}$ and $1 \leq i \leq n$,  denote the $(i,i)^{th}$-entry of $\lambda$ by $\lambda_i$.   Given $\lambda_1,\lambda_2,\cdots,\lambda_n \in \mathbb{T}$, the diagonal matrix with diagonal entries $\lambda_1,\lambda_2,\cdots,\lambda_n$ will be denoted by $diag(\lambda_1,\lambda_2,\cdots,\lambda_n)$. Here  by $\mathbb{T}$, we mean the unit circle of the complex plane.  
 
 For $\sigma \in S_{n}$, let $U_{\sigma}$ be the unitary on $\mathbb{C}^{n}$ such that $U_{\sigma}(e_i)=e_{\sigma(i)}$ for $i \in \{1,2,\cdots,n\}$. For $\lambda:=diag(\lambda_1,\lambda_2,\cdots,\lambda_n) \in \mathbb{T}^{n}$ and $\sigma \in S_n$, let $\lambda_{\sigma}:=diag(\lambda_{\sigma(1)},\lambda_{\sigma(2)},\cdots,\lambda_{\sigma(n)})$. Note that for $\sigma \in S_n$ and $\lambda \in \mathbb{T}^{n}$, $U_{\sigma}\lambda U_{\sigma}^{*}=\lambda_{\sigma^{-1}}$. 
 The following proposition may be  known to experts.   We use the notation introduced in the preceding two paragraphs in the proof of the following proposition.
 \begin{ppsn}
 \label{representation of the unitary group}
Let $n \geq 2$,  $\clh$ be a Hilbert space and $\rho$ be a strongly continuous unitary representation of $U(n)$ on $\clh$. Suppose that the dimension of $\clh$ is strictly less than $n$. Then $\rho(U)=1$ for every $U \in SU(n)$. Here $1$ denotes the identity operator on $\clh$. 
  \end{ppsn}
 \textit{Proof.} For $m:=(m_1,m_2,\cdots,m_n) \in \mathbb{Z}^{n}$, let \[\clh_{m}:=\{v \in \clh: \textit{~for every}~\lambda \in \mathbb{T}^{n}, ~\rho(\lambda)v=\lambda_{1}^{m_1}\lambda_{2}^{m_2}\cdots \lambda_{n}^{m_n}v\}.\] Restrict $\rho$ to the compact abelian group $\mathbb{T}^{n}$. Then the Hilbert space $\clh$ decomposes as $\displaystyle \clh=\bigoplus_{m \in \mathbb{Z}^{n}}\clh_{m}$. 
 
 Fix $m:=(m_1,m_2,\cdots,m_n) \in \mathbb{Z}^{n}$. Suppose $\clh_m \neq 0$. Then $m_i=m_j$ for every $i,j \in \{1,2,\cdots,n\}$. Suppose not. Then there exists $i,j \in \{1,2,\cdots,n\}$ such that $i \neq j$ and $m_i \neq m_j$. Let $v \in \clh_m$,  $\lambda:=diag(\lambda_1,\lambda_2,\cdots,\lambda_n) \in \mathbb{T}^{n}$ and $\sigma \in S_n$ be given. Calculate as follows to observe that 
 \begin{align*}
 \rho(U_{\sigma})\rho(\lambda)\rho(U_{\sigma})^{*}v&=\rho(U_{\sigma}\lambda U_{\sigma}^{*})v \\
                                                                                  &=\rho(\lambda_{\sigma^{-1}})v \\
                                                                                  &=\lambda_{\sigma^{-1}(1)}^{m_1}\lambda_{\sigma^{-1}(2)}^{m_2}\cdots\lambda_{\sigma^{-1}(n)}^{m_n}v \\ 
                                                                                  &=\lambda_{1}^{m_{\sigma(1)}}\lambda_{2}^{m_{\sigma(2)}}\cdots \lambda_{n}^{m_{\sigma(n)}}v.
  \end{align*}
  The  calculation implies that $\rho(U_{\sigma})^{*}$ maps $\clh_m$ into $\clh_{m_{\sigma}}$. This implies in particular that $\clh_{m_\sigma} \neq 0$ for every $\sigma \in S_n$. Note that for $m^{'},m^{''} \in \mathbb{Z}^{n}$ if $m^{'} \neq m^{''}$ then $\clh_{m^{'}}$ is orthogonal to $\clh_{m^{''}}$. This together with the fact that the cardinality of  $\{m_{\sigma}: \sigma \in S_n\}$ is at least $n$ (Lemma \ref{combinatorial lemma}) implies that the dimension of $\clh$ is at least $n$ which contradicts the hypothesis. Hence if $\clh_m \neq 0$ then $m_i=m_j$ for every $i,j \in \{1,2,\cdots,n\}$. This has the consequence that if $\lambda \in \mathbb{T}^{n} \cap SU(n)$ then $\rho(\lambda)=1$. 
  
  Let $U \in SU(n)$ be given. Then there exists $\lambda \in \mathbb{T}^{n}$ and $V \in U(n)$ such that $V\lambda V^{*}=U$. Since $U \in SU(n)$, it follows that the determinant of $\lambda$ is one. Hence \[
  \rho(U)=\rho(V\lambda V^{*})=\rho(V)\rho(\lambda)\rho(V)^{*}=\rho(V)\rho(V)^{*}=\rho(VV^{*})=\rho(1)=1.\]
  This completes the proof. \hfill $\Box$
  
  Let $G$ be a compact group, let $\clh$ be a separable Hilbert space and let $\pi:G \to B(\clh)$ be a strongly continuous unitary representation of $G$ on $\clh$. Let $\{(\clh_{\alpha},\pi_{\alpha})\}_{\alpha \in \Lambda}$ be the complete list of irreducible subrepresentations occuring in $(\clh,\pi)$. For $\alpha \in \Lambda$, $\clh_{\alpha}$ is finite dimensional. For $\alpha \in \Lambda$, let $n_{\alpha}$ be the multiplicity of $(\clh_{\alpha},\pi_{\alpha})$ in $(\clh,\pi)$. For $\alpha \in \Lambda$, let $\ell_{2}^{n_{\alpha}}$ be a Hilbert space of dimension $n_{\alpha}$. Up to a unitary equivalence, we can write 
  \begin{align*}
  \clh&=\bigoplus_{\alpha \in \Lambda}(\clh_{\alpha} \otimes \ell_{2}^{n_\alpha}) \\
  \pi(g)&=\bigoplus_{\alpha \in \Lambda}(\pi_{\alpha}(g) \otimes 1)
    \end{align*}
  for $g \in G$. Let \[\displaystyle \bigoplus_{\alpha \in \Lambda}B(\clh_{\alpha}):=\big\{(T_{\alpha})_{\alpha \in \Lambda}: \sup_{\alpha \in \Lambda}||T_{\alpha}||<\infty\big\}.\] For $\displaystyle T:=(T_{\alpha})_{\alpha \in \Lambda} \in \bigoplus_{\alpha \in \Lambda}B(\clh_{\alpha})$, define \[||T||:=\displaystyle \sup_{\alpha \in \Lambda}||T_{\alpha}||.\] Then $\Big(\displaystyle \bigoplus_{\alpha \in \Lambda}B(\clh_{\alpha}),||~||\Big)$ is a $C^{*}$-algebra. For $T=(T_\alpha)_{\alpha \in \Lambda} \in \displaystyle \bigoplus_{\alpha \in \Lambda}B(\clh_\alpha)$, define $\widetilde{T} \in B(\clh)$ by the formula \[\widetilde{T}:=\displaystyle \bigoplus_{\alpha \in \Lambda}(T_{\alpha} \otimes 1).\] The map $\displaystyle \bigoplus_{\alpha \in \Lambda}B(\clh_{\alpha}) \ni T \to \widetilde{T} \in B(\clh)$ is an injective $*$-homomorphism. The proof of the following proposition is elementary. Thus we omit its proof.
  
  \begin{ppsn}
  \label{von Neumann algebra generated by a compact group}
  With the foregoing notation, we have 
  \[
  \pi(G)^{''}=\{\widetilde{T}: T \in \bigoplus_{\alpha \in \Lambda}B(\clh_\alpha)\}.\]
    \end{ppsn}

 Let $\clh_1,\clh_2,\clk_1,\clk_2$ be non-zero separable Hilbert spaces. For $i \in \{1,2\}$, let $M_i$ be a unital commutative von Neumann algebra acting on $\clh_i$. Fix $i \in \{1,2\}$. Consider the tensor product von Neumann algebra $M_i \otimes B(\clk_i)$ acting on $\clh_i \otimes \clk_i$. Denote the unitary group of $M_i \otimes B(\clk_i)$ endowed with the strong operator topology by $\mathcal{U}(M_i \otimes B(\clk_i))$. Suppose that $dim(\clk_1)<dim(\clk_2)$. Let $\widetilde{\clk_2}$ be a finite dimensional subspace of $\clk_2$ such that $dim(\clk_1)<dim(\widetilde{\clk_2})$. Denote the unitary group of $\widetilde{\clk_2}$ by $\mathcal{U}(\widetilde{\clk_2})$ and endow $\mathcal{U}(\widetilde{\clk_2})$ with the norm topology. Write $\clk_2=\widetilde{\clk_2}\oplus \widetilde{\clk_2}^{\perp}$. We denote the special unitary group of $\widetilde{\clk_2}$ by $SU(\widetilde{\clk_2})$, i.e. 
 \[
 SU(\widetilde{\clk_2}):=\{U \in \mathcal{U}(\widetilde{\clk_2}):~\det(U)=1\}.\] For $U \in \mathcal{U}(\widetilde{\clk_2})$, define $\widetilde{U} \in B(\clk_2)$ by $\widetilde{U}=U \oplus 1$. Observe  that the map \[\mathcal{U}(\widetilde{\clk_2}) \ni U \to 1 \otimes \widetilde{U} \in \mathcal{U}(M_{2} \otimes B(\clk_2))\] is a topological embedding and is also a group homomorphism. We use the preceeding notation in the statement and the proof of the following proposition.
 
 \begin{ppsn}
 \label{Sunder}
Let $\Phi:\mathbb{R}^{d} \times \mathcal{U}(M_2 \otimes B(\clk_2)) \to \mathbb{R}^{d} \times \mathcal{U}(M_1 \otimes B(\clk_1))$ be a continuous group homomorphism. Then $\{(0,1\otimes \widetilde{U}): U \in SU(\widetilde{\clk_2})\}$ is contained in the kernel of $\Phi$. In particular, $\Phi$ is not $1$-$1$. 

Here, for $i\in \{1,2\}$, the topology on  $\mathbb{R}^{d} \times \mathcal{U}(M_i \otimes B(\clk_i))$ is the product topology and the group structure on $\mathbb{R}^{d} \times \mathcal{U}(M_i \otimes B(\clk_i))$ is that  of cartesian product. 
 
 \end{ppsn}
 \textit{Proof.} Let $\pi_1:\mathbb{R}^{d} \times \mathcal{U}(M_1 \otimes B(\clk_1)) \to \mathbb{R}^{d}$ and $\pi_{2}:\mathbb{R}^{d} \times \mathcal{U}(M_1 \otimes B(\clk_1)) \to \mathcal{U}(M_1 \otimes B(\clk_1))$ be the first and second co-ordinate projections. Define $\Phi_1=\pi_1 \circ \Phi$ and $\Phi_2=\pi_2 \circ \Phi$. 
  Let $\mathcal{L}:=\clh_1 \otimes \clk_1$ and $\pi:\mathcal{U}(\widetilde{\clk_2})\to B(\mathcal{L})$ be the strongly continuous unitary representation defined by the equation
 \[
 \pi(U):=\Phi_{2}(0, 1\otimes \widetilde{U}).\]
  Let $\{(\mathcal{L}_{\alpha},\pi_{\alpha})\}_{\alpha \in \Lambda}$ be the complete list of irreducible subrepresentations of $\mathcal{U}(\widetilde{\clk_2})$ occuring in $(\mathcal{L},\pi)$. Note that for each $\alpha \in \Lambda$, $\mathcal{L}_{\alpha}$ is finite dimensional. For $\alpha \in \Lambda$, let $n_{\alpha}$ be the multiplicity of $(\mathcal{L}_{\alpha},\pi_{\alpha})$ in $(\mathcal{L},\pi)$. We use/apply the notation explained before Prop. \ref{von Neumann algebra generated by a compact group} to the compact group $\mathcal{U}(\widetilde{\clk_2})$ and the representation $(\mathcal{L},\pi)$. 
 
 \textit{Claim:} $dim(\mathcal{L}_{\alpha}) \leq dim(\mathcal{K}_{1})$ for every $\alpha \in \Lambda$. Let $\alpha_0 \in \Lambda$ be fixed. For $S \in B(\mathcal{L}_{\alpha_0})$, let $\displaystyle \underline{S}:=(\underline{S}_{\alpha})_{\alpha \in \Lambda} \in \bigoplus_{\alpha \in \Lambda}B(\mathcal{L}_{\alpha})$ be defined by the following equation
 \begin{equation*}
\underline{S}_{\alpha}:=\begin{cases}
 S  & \mbox{ if
} \alpha=\alpha_0,\cr
   &\cr
   0 &  \mbox{ if } \alpha \neq \alpha_0.
         \end{cases}
\end{equation*}
Clearly the map $\displaystyle B(\mathcal{L}_{\alpha_0}) \ni S \to \underline{S} \in \bigoplus_{\alpha \in \Lambda} B(\mathcal{L}_{\alpha})$ is an injective $*$-homomorphism. Let $P \in B(\mathcal{L}_{\alpha_0})$ be a non-zero element.

 Since $\pi(\mathcal{U}(\widetilde{\clk_2})) \subset \mathcal{U}(M_1\otimes B(\clk_1))$, it follows that the von Neumann algebra generated by $\pi(\mathcal{U}(\widetilde{\clk_2}))$, i.e. $\pi(\mathcal{U}(\widetilde{\clk_2})) ^{''}$, is contained in $M_1 \otimes B(\clk_1)$. 
 Note that  $\widetilde{\underline{P}} \in M_1 \otimes B(\clk_1)$ is non-zero. Treating $M_1$ as a $C^{*}$-algebra, we see that there exists a character $\chi$ of $M_1$ such that $(\chi \otimes 1)(\widetilde{\underline{P}}) \neq 0$. This shows that map $B(\mathcal{L}_{\alpha_0}) \ni S \to (\chi \otimes 1)(\widetilde{\underline{S}}) \in B(\clk_1)$ is a non-zero $*$-homomorphism. Since $\mathcal{L}_{\alpha_0}$ is finite dimensional, it follows that $B(\mathcal{L}_{\alpha_0})$ has no non-zero ideal. As a consequence, we conclude that the map $B(\mathcal{L}_{\alpha_0}) \ni S \to (\chi \otimes 1)(\widetilde{\underline{S}}) \in B(\clk_1)$ is an injection. Hence $dim(\mathcal{L}_{\alpha_0}) \leq dim(\clk_1)$. This proves our claim. 
 
 Thanks to Prop. \ref{representation of the unitary group} and to the fact that $dim(\mathcal{L}_{\alpha})<dim(\widetilde{\clk_2})$ for every $\alpha \in \Lambda$, we obtain that  $\pi(U)=1$ for every $U \in SU(\widetilde{\clk_2})$. This implies that $\Phi_2(0,1\otimes \widetilde{U})=1$ for every $U \in SU(\widetilde{\clk_2})$.  Note that the map $\Phi_1:\mathbb{R}^{d} \times \mathcal{U}(M_{2} \otimes B(\clk_2)) \to \mathbb{R}^{d}$ is continuous. Consequently $\{\Phi_1(0,1 \otimes \widetilde{U}): U \in SU(\widetilde{\clk_2})\}$ is a compact subgroup of $\mathbb{R}^{d}$. But the only compact subgroup of $\mathbb{R}^{d}$ is the trivial one i.e. $\{0\}$. Hence $\Phi_1(0,1 \otimes \widetilde{U})=0$ for every $U \in SU(\widetilde{\clk_2})$.
 As a consequence, we obtain that $\{(0,1\otimes \widetilde{U}): U \in SU(\widetilde{\clk_2})\}$ is contained in the kernel of $\Phi$. This completes the proof. \hfill $\Box$
 
 The following corollary is immediate from Prop. \ref{Sunder} and Corollary \ref{commutant 2}.
 
 \begin{crlre}
 Let $A$ be a $P$-module and $k_1,k_2 \in \{1,2,\cdots,\}\cup \{\infty\}$. For $i \in \{1,2\}$, denote the CCR flow associated to the $P$-module $A$ of multiplicity $k_i$ by $\alpha^{(A,k_i)}$. Then $\alpha^{(A,k_1)}$ is cocycle conjugate to $\alpha^{(A,k_2)}$ if and only if $k_1=k_2$.
  \end{crlre}
  
   Let $\clh_1,\clh_2,\clk$ be non-zero separable Hilbert spaces. Assume that $\clk$ is infinite dimensional. For $i \in \{1,2\}$, let $M_i \subset B(\clh_i)$ be a unital commutative von Neumann algebra. For $i \in \{1,2\}$, consider the tensor product von Neumann algebra $M_i \otimes B(\clk)$ acting on $\clh_i \otimes \clk$. The unitary groups of $M_1, M_1 \otimes B(\clk)$ and $M_2 \otimes B(\clk)$ are endowed with the corresponding strong operator topologies  and we will denote them by $\mathcal{U}(M_1)$, $\mathcal{U}(M_1 \otimes B(\clk))$ and $\mathcal{U}(M_2 \otimes B(\clk))$ respectively. We denote the identity element of various unitary groups involved by $1$. The identity element of $\mathbb{R}^{d}$ will be denoted by $0$.
  
  \begin{lmma}
  \label{Sundar}
  With the foregoing notation, the topological groups $\mathbb{R}^{d} \times \mathcal{U}(M_1 \otimes B(\clk))$ and $\mathbb{R}^{d} \times \mathcal{U}(M_1) \times \mathcal{U}(M_2 \otimes B(\clk))$ are not isomorphic. 
  \end{lmma}
  \textit{Proof.} Suppose that there exists a map, say, \[\Phi: \mathbb{R}^{d} \times \mathcal{U}(M_1 \otimes B(\clk)) \to \mathbb{R}^{d} \times \mathcal{U}(M_1) \times \mathcal{U}(M_2 \otimes B(\clk))\] such that $\Phi$ is a topological group isomorphism. Denote the second co-ordinate projection from $\mathbb{R}^{d} \times \mathcal{U}(M_1) \times \mathcal{U}(M_2 \otimes B(\clk))$ onto $\mathcal{U}(M_1)$ by $\pi$. 
  Define $\widetilde{\Phi}:\mathcal{U}(M_1 \otimes B(\clk))\to \mathcal{U}(M_1)$ as follows: for $U \in \mathcal{U}(M_1 \otimes B(\clk))$, let $\widetilde{\Phi}(U)=\pi \circ \Phi(0,U)$.  Note that $\widetilde{\Phi}$ is a continuous group homomorphism. 
  
  We claim that for $x \in \mathcal{U}(M_1)$, $\widetilde{\Phi}(x\otimes 1)=1$. Let $x \in \mathcal{U}(M_1)$ be given. 
    Choose an orthonormal basis, say,  $\{\xi_1,\xi_2,\xi_3, \cdots\}$ of $\clk$. Let $\mathbb{N}:=\{1,2,3,\cdots\}$. For $n \in \mathbb{N}$, let $E_{n}$ be the orthogonal projection onto the $1$-dimensional subspace of $\clk$ spanned by $\{\xi_n\}$. For $m,n \in \mathbb{N}$, let $U_{m,n}$ be a unitary on $\clk$ such that $U_{m,n}E_nU_{m,n}^{*}=E_m$. For $n \in \mathbb{N}$, define 
  \begin{align*}
  T_{n}:&=x\otimes E_n+1 \otimes (1-E_n) \\
  S_{n}:&=T_1T_2\cdots T_n
    \end{align*}
    Note that for $n \in \mathbb{N}$, $S_n=x \otimes (\sum_{k=1}^{n}E_k)+1 \otimes (1-\sum_{k=1}^{n}E_{k})$. Note that the sequence $\{T_{n}\}_{n \in \mathbb{N}}$ converges strongly to $1$ and the sequence $\{S_n\}_{n \in \mathbb{N}}$ converges strongly to $x \otimes 1$. 
    
    Fix $m,n \in \mathbb{N}$. Note that $(1\otimes U_{m,n})T_{n}(1 \otimes U_{m,n})^{-1}=T_m$. The fact that $\mathcal{U}(M_1)$ is abelian implies that $\widetilde{\Phi}(T_n)=\widetilde{\Phi}(T_m)$. Hence the sequence $\{\widetilde{\Phi}(T_n)\}_{n \in \mathbb{N}}$ is a constant sequence.  
    Since $\{T_n\} \to 1$ and $\widetilde{\Phi}$ is continuous, it follows that $\widetilde{\Phi}(T_n)=1$ for every $n \in \mathbb{N}$. The fact that $\widetilde{\Phi}$ is a group homomorphism implies that $\widetilde{\Phi}(S_n)=1$ for every $n \in \mathbb{N}$. But $\{S_n\} \to x \otimes 1$ and $\widetilde{\Phi}$ is continuous. As a consequence, we obtain that $\widetilde{\Phi}(x \otimes 1)=1$. This proves our claim. 
    
    Since $\Phi$ is a group isomorphism, it follows that $\Phi$ maps the center of the topological group $\mathbb{R}^{d} \times \mathcal{U}(M_1 \otimes B(\clk))$, which is $\{(\lambda,x \otimes 1): \lambda \in \mathbb{R}^{d}, x \in \mathcal{U}(M_1)\}$, onto the center of $\mathbb{R}^{d} \times \mathcal{U}(M_1) \times \mathcal{U}(M_2 \otimes B(\clk))$, which is $\{(\mu,y,z \otimes 1):\mu \in \mathbb{R}^{d}, y \in \mathcal{U}(M_1), z \in \mathcal{U}(M_2)\}$. 
    Consider the element $(0,-1,1\otimes 1) \in \mathbb{R}^{d} \times \mathcal{U}(M_1) \times \mathcal{U}(M_2 \otimes B(\clk))$ which is in the center of $\mathbb{R}^{d} \times \mathcal{U}(M_1) \times \mathcal{U}(M_2 \otimes B(\clk))$.  Hence there exists $\lambda \in \mathbb{R}^{d}$ and $x \in \mathcal{U}(M_1)$ such that $\Phi(\lambda,x\otimes 1)=(0,-1,1\otimes 1)$. Note that $(0,-1,1\otimes 1)$ has order $2$. Since $\Phi$ is a group isomorphism, it follows that $(\lambda,x\otimes 1)$ has order $2$. This implies that $\lambda=0$. Consequently, we have $\widetilde{\Phi}(x \otimes 1)=-1$ which is a contradiction to the fact that $\widetilde{\Phi}(x\otimes 1)=1$. Hence the proof. \hfill $\Box$
  
  Let $A_1,A_2$ be $P$-modules and $k_1,k_2 \in \{1,2,\cdots,\}\cup\{\infty\}$. Fix $i \in \{1,2\}$. Let $\clk_i$ be a Hilbert space of dimension $k_i$. Let $V^{(i)}:P \to B(L^{2}(A_i)\otimes \clk_i)$ be the isometric representation associated to the $P$-module $A_i$ of multiplicity $k_i$.  Denote the isometric representation associated to the $P$-module $A_i$ of multiplicity $1$ by $\widetilde{V^{(i)}}$. Set $\clh_i:=L^{2}(A_i)\otimes \clk_i$, $\clh:=\clh_1 \oplus \clh_2$ and $V:=V^{(1)}\oplus V^{(2)}$. Let $\{W^{(i)}_{z}\}_{z \in \mathbb{R}^{d}}$ be the family of partial isometries, described in the paragraph following Theorem \ref{main theorem of the appendix}, associated to the isometric representation $V^{(i)}$, and let $\{\widetilde{W^{(i)}_{z}}\}_{z \in \mathbb{R}^{d}}$ and $\{W_{z}\}_{z \in \mathbb{R}^{d}}$ be the family of partial isometries  associated to the isometric representations $\widetilde{V^{(i)}}$ and $V$ respectively. Note that $W^{(i)}_{z}=\widetilde{W^{(i)}_{z}}\otimes 1$. We use the notation developed in the paragraphs between Lemma \ref{von Neumann algebra generated by Wiener-Hopf operators} and Proposition \ref{Induced representation}. 
  
  Let $G_{A_i}$ be the isotropy group of $A_i$, i.e. $G_{A_i}:=\{z \in \mathbb{R}^{d}: A_{i}+z=A_i\}$. For $z \in G_{A_i}$, let $U^{(i)}_{z}$ be the unitary on $L^{2}(A_i)$ defined by the formula
  \[
  U^{(i)}_{z}f(x):=f(x-z)
  \]
  for $f \in L^{2}(A_i)$. 
  Let $M_i$ be the von Neumann algebra generated by $\{U^{(i)}_{z}:z \in G_{A_i}\}$ acting on $L^{2}(A_i)$. Denote the CCR flow associated to the $P$-module $A_i$ of multiplicity $k_i$ by $\alpha^{(A_i,k_i)}$.   
 \begin{ppsn}
 \label{Gauge group of the tensor product}
 Assume that $A_2$ is not a translate of $A_1$, i.e. for every $z \in \mathbb{R}^{d}$, $A_{1}+z \neq A_2$.
 With the foregoing notation, the gauge group of  $\alpha^{(A_1,k_1)}\otimes \alpha^{(A_2,k_2)}$  is isomorphic to  $\mathbb{R}^{d} \times \mathcal{U}(M_1 \otimes B(\clk_1)) \times \mathcal{U}(M_2 \otimes B(\clk_2))$. 
 \end{ppsn}
 \textit{Proof.} Note that $\alpha^{(A_1,k_1)}\otimes \alpha^{(A_2,k_2)}$ is  the CCR flow, denoted $\alpha^{V}$,  associated to the isometric representation $V$. Since $V^{(1)}$ and $V^{(2)}$ admits no non-zero additive cocycles, by Remark \ref{additive cocycles}, it follows that $V$ admits no non-zero additive cocycle. By Theorem \ref{description of the gauge group}, it follows that the gauge group of $\alpha^{V}$ is isomorphic to $\mathbb{R}^{d} \times \mathcal{U}(M)$ where $M$ is the commutant of the von Neumann algebra generated by $\{V_{x}:x \in P\}$. 
 
 By Corollary \ref{commutant 2}, for $i \in \{1,2\}$, the commutant of the von Neumann algebra generated by $\{V_{x}^{(i)}:x \in P\}$ is $M_i \otimes B(\clk_i)$. We write operators acting on $\clh=\clh_1\oplus \clh_2$ in terms of block matrices. We claim that 
 \[
 M=\Big\{ \begin{pmatrix}
                 T_1 & 0 \\
                 0 & T_2 \\
                 \end{pmatrix}: T_1 \in M_{1} \otimes B(\clk_1), T_2 \in M_{2} \otimes B(\clk_2)\Big\}
  \]
  Once the above claim is established,  thanks  to Theorem \ref{description of the gauge group}, the conclusion follows immediately.
 It is clear that $
 \Big\{ \begin{pmatrix}
                 T_1 & 0 \\
                 0 & T_2 \\
                 \end{pmatrix}: T_1 \in M_{1}\otimes B(\clk_1), T_2 \in M_{2}\otimes B(\clk_2)\Big\}$ is contained in $M$. Let $T:=\begin{pmatrix}
                                              T_{11} & T_{12} \\
                                              T_{21} & T_{22}
                                              \end{pmatrix} \in M$ be given. It is routine to verify that $T_{12}W^{(2)}_{z}=W^{(1)}_{z}T_{12}$ for $z \in \mathbb{R}^{d}$. Let $f \in C_{c}(\mathbb{R}^{d})$ be given. Calculate as follows to observe that 
  \begin{align*}
  T_{12}(\pi_{A_2}(\widetilde{f})\otimes 1)&=T_{12}\Big(\Big(\int f(-z)\widetilde{W^{(2)}_{z}}dz\Big)\otimes 1 \Big)\\
                                                                 &=T_{12}\Big(\int f(-z)W^{(2)}_{z}dz\Big)\\
                                                                 &=\int f(-z)T_{12}W^{(2)}_{z}dz\\
                                                                 &=\int f(-z)W^{(1)}_{z}T_{12}dz\\
                                                                 &=\Big(\int f(-z)W^{(1)}_{z}dz\Big)T_{12}\\
                                                                 &=\Big(\int f(-z)(\widetilde{W^{(1)}_{z}} \otimes 1)dz \Big)T_{12}\\
                                                                 &=(\pi_{A_1}(\widetilde{f})\otimes 1)T_{12}.
    \end{align*}
    Since $\{\widetilde{f}:f \in C_{c}(\mathbb{R}^{d})\}$ generates $C^{*}(\clg_u)$, it follows that $T_{12}$ intertwines the representation $(\pi_{A_2}(.)\otimes 1, \clh_2)$ and the representation $(\pi_{A_1}(.)\otimes 1,\clh_1)$. Theorem \ref{main theorem of the appendix} and  the hypothesis  $A_1+z \neq A_2$ for every $z \in \mathbb{R}^{d}$ implies that $\pi_{A_1}$ and $\pi_{A_2}$ are disjoint. This implies that $\pi_{A_1}(.)\otimes 1$ and $\pi_{A_2}(.)\otimes 1$ are disjoint. Consequently $T_{12}=0$. In a similar fashion, we conclude that $T_{21}=0$. It is now clear that $T_{11} \in M_{1}\otimes B(\clk_1)$ and $T_{22} \in M_{2}\otimes B(\clk_2)$. This proves our claim. Hence the proof. \hfill $\Box$
                                              
  Now we prove  Theorem \ref{main theorem}. We use the notation developed in the two paragraphs that precede Proposition \ref{Gauge group of the tensor product} and we write $\alpha \cong \beta$ to indicate that $\alpha$ and $\beta$ are cocycle conjugate.
  
  \textit{Proof of Theorem \ref{main theorem}.} As mentioned in the introduction, it is clear that $(3) \implies (1) \implies (2)$. Suppose that $(2)$ holds. Then the gauge group of $\alpha^{(A_1,k_1)}$ is isomorphic to the gauge group of $\alpha^{(A_2,k_2)}$. By Corollary \ref{commutant 2}, it follows that $\mathbb{R}^{d} \times \mathcal{U}(M_1 \otimes B(\clk_1))$ is isomorphic to $\mathbb{R}^{d} \times \mathcal{U}(M_2 \otimes B(\clk_2))$. Since $M_1$ and $M_2$ are abelian, it follows from Proposition \ref{Sunder} that $k_1=k_2$. With no loss of generality, we can assume that $\clk_1=\clk_2$. 
 
 Suppose, on the contrary, assume that $A_1$ and $A_2$ are not translates of each other. Since $\alpha^{(A_1,k_1)}$ is cocycle conjugate to $\alpha^{(A_2,k_2)}$, it follows that \[\alpha^{(A_1,k_1+1)}\cong \alpha^{(A_1,1)} \otimes \alpha^{(A_1,k_1)}\cong \alpha^{(A_1,1)}\otimes \alpha^{(A_2,k_2)}.\] Hence $\alpha^{(A_1,k_1+1)}$ and $\alpha^{(A_1,1)}\otimes \alpha^{(A_2,k_2)}$ have isomorphic gauge groups. Corollary \ref{commutant 2} and Proposition \ref{Gauge group of the tensor product} together imply  that the topological groups $\mathbb{R}^{d} \times \mathcal{U}(M_1 \otimes B(\widetilde{\clk_1}))$ and $\mathbb{R}^{d} \times \mathcal{U}(M_1) \times \mathcal{U}(M_2 \otimes B(\clk_2))$ are isomorphic, where $\widetilde{\clk_1}$ is a Hilbert space of dimension $k_1+1$. Let $\Phi: \mathbb{R}^{d} \times \mathcal{U}(M_1 \otimes B(\widetilde{\clk_1}))\to \mathbb{R}^{d} \times \mathcal{U}(M_1) \times \mathcal{U}(M_2 \otimes B(\clk_2))$ be a topological group isomorphism. 
 
 Suppose $k_1=k_2$ is finite. Let $\pi_{12}:\mathbb{R}^{d} \times \mathcal{U}(M_1) \times \mathcal{U}(M_2 \otimes B(\clk_2))\to \mathbb{R}^{d} \times \mathcal{U}(M_1) $ and $\pi_{13}:\mathbb{R}^{d} \times \mathcal{U}(M_1) \times \mathcal{U}(M_2 \otimes B(\clk_2)) \to \mathbb{R}^{d} \times \mathcal{U}(M_2 \otimes B(\clk_2))$ be defined by the following formulas
 \begin{align*}
 \pi_{12}(x,Y,Z)&=(x,Y) \\
 \pi_{13}(x,Y,Z)&=(x,Z)
  \end{align*}
for $(x,Y,Z) \in \mathbb{R}^{d} \times \mathcal{U}(M_1) \times \mathcal{U}(M_2 \otimes B(\clk_2))$. Let $\Phi_{12}=\pi_{12} \circ \Phi$ and $\Phi_{13}=\pi_{13} \circ \Phi$. Proposition \ref{Sunder} implies that $\{(0,1\otimes U): U \in SU(\widetilde{\clk_1})\}$ is contained in the kernel of both $\Phi_{12}$ and $\Phi_{23}$. This implies that $\{(0,1\otimes U): U \in SU(\widetilde{\clk_1})\}$ is contained in the kernel of $\Phi$, which is a contradiction. This implies that $k_1=k_2=\infty$. Hence $\widetilde{\clk_1}$ and $\clk_2$ are infinite dimensional separable Hilbert spaces. The fact that $\Phi$ is a topological group isomorphism is a contradiction to Lemma \ref{Sundar}. These contradictions are due to our initial assumption that $A_1$ and $A_2$ are not translates of each other. Hence there exists $z \in \mathbb{R}^{d}$ such that $A_1+z=A_2$. The proof of the implication $(2) \implies (3)$ is now complete. \hfill $\Box$
 
 \section{Appendix}
 
 Here we provide a proof of  Theorem \ref{main theorem of the appendix}. The proof is an application of Rieffel's theory of Morita equivalence. Rieffel's theorem, Theorem 6.23 of \cite{Rieffel}, asserts that if $A$ and $B$ are Morita equivalent C*-algebras then the 
 category of representations of $A$ and that of $B$ are equivalent. As the Deaconu-Renault groupoid $X \rtimes P$ is equivalent to a transformation groupoid $Y \rtimes \mathbb{R}^{d}$, it follows from \cite{MRW} that the $C^{*}$-algebras $C^{*}(X \rtimes P)$ and $C^{*}(Y \rtimes \mathbb{R}^{d}) \cong C_{0}(Y) \rtimes \mathbb{R}^{d}$ are Morita equivalent. Consequently it suffices to prove the result for a transformation groupoid.
 
 Let us begin by reviewing the basics of Rieffel's notion of Morita equivalence. Let $B$ be a $C^{*}$-algebra. For a Hilbert $B$-module $E$, we denote the $C^{*}$-algebra of adjointable operators on $E$ by $\mathcal{L}_{B}(E)$ and the $C^{*}$-algebra of compact operators by $\mathcal{K}_{B}(E)$. The Hilbert module $E$ is said to be full if the linear span of $\{\langle x|y \rangle: x,y \in E\}$ is dense in $B$. 
 Let $A$ and $B$ be $C^{*}$-algebras. By an $A$-$B$ \emph{imprimitivity bimodule}, we mean a pair $(E,\phi)$, where $E$ is a full Hilbert $B$-module, $\phi:A \to \mathcal{L}_{B}(E)$ is an injective $*$-homomorphism and $\phi(A)=\mathcal{K}_{B}(E)$. 
 The $C^{*}$-algebras $A$ and $B$ are said to be Morita equivalent if there exists an $A$-$B$ imprimitivity bimodule.
 
  Next we recall Rieffel's \emph{induction} procedure. Let $A$ and $B$ be Morita equivalent $C^{*}$-algebras with $E$ being an $A$-$B$ imprimitivity bimodule. Suppose $\pi$ is a representation of $B$ on a Hilbert space $\mathcal{H}_{\pi}$. Consider the internal tensor product $E \otimes_{\pi} \clh_{\pi}$ which is.a Hilbert space. The $C^{*}$-algebra $A$ acts on the Hilbert space $E \otimes_{\pi} \clh_{\pi}$ as follows: for $a \in A$, let $Ind(\pi)(a):=\phi(a) \otimes 1$. Then $Ind(\pi)$ is a representation of $A$. Rieffel's fundamental theorem, Theorem 6.23 of \cite{Rieffel}, asserts that
  \[
  \pi \to Ind(\pi)
  \]
  is a functor which identifies the category of representations of $B$ and the category of representations of $A$. Let us isolate two consequences of the above fact in the following remark.
  
  \begin{rmrk}
  \label{Rieffel commutant}
  With the foregoing notation, we have the following.
  \begin{enumerate}
  \item[(1)] Let $\pi_1$ and $\pi_2$ be representations of $B$. Then $\pi_1$ and $\pi_2$ are disjoint if and only if $Ind(\pi_1)$ and $Ind(\pi_2)$ are disjoint.
  \item[(2)] Let $\pi$ be a representation of $B$. Then $(Ind(\pi)(A))^{'}=\{1 \otimes F: F \in \pi(B)^{'}\}$.
  
  \end{enumerate}

  \end{rmrk}
  
 We keep the notation explained in the paragraphs starting from line 8, Page 17 until the end of Theorem \ref{main theorem of the appendix}.
  Let us fix a few notation:  Let $Z:=\{(y,s):y+s \in X\}$ and let $\rho:Z \to Y$ and $\sigma:Z \to X$ be defined by $\rho(y,s)=y$ and $\sigma(y,s)=y+s$. Then $Z$ is a $(Y \rtimes \mathbb{R}^{d},X \rtimes P)$-equivalence where the actions of $Y \rtimes \mathbb{R}^{d}$ and $X \rtimes P$ on $Z$ are given by the formulas: 
  \begin{align*}
(y,s)(z,t)&=(y,s+t) \textit{~if~~$y+s=z$~}\\
(z,t)(x,r)&=(z,t+r) \textit{~if ~~$z+t=x$~}
\end{align*}  
for $(x,r) \in X \rtimes P$, $(y,s) \in Y \rtimes \mathbb{R}^{d}$ and $(z,t) \in Z$. For the definition of an action of a groupoid on a space and for the notion of groupoid equivalence, we refer the reader to \cite{MRW}.  Let $A:=C^{*}(Y \rtimes \mathbb{R}^{d})$ and $B:=C^{*}(X \rtimes P)$. Denote $C_{c}(Y \rtimes \mathbb{R}^{d})$ and $C_{c}(X \rtimes P)$ by $\mathcal{A}$ and $\mathcal{B}$ respectively. Note that $\mathcal{A}$ and $\mathcal{B}$ are dense in $A$ and $B$ respectively.  Denote $C_{c}(Z)$ by $\mathcal{E}$. For $\xi \in \mathcal{A}$, $\chi \in \mathcal{E}$, and $\eta \in \mathcal{B}$, let 
\begin{align*}
(\xi.\chi)(y,r)&=\int \xi(y,s)\chi(y+s,r-s)ds, \\
(\chi.\eta)(y,r)&=\int \chi(y,s)\eta(y+s,r-s)1_{X}(y+s)ds, \textrm{~and~} \\
\langle \chi_1,\chi_2\rangle_{\mathcal{B}}(x,r)&=\int \overline{\chi_1(x+s,-s)}\chi_2(x+s,r-s)ds.
\end{align*}
The above formulas make $\mathcal{E}$ into a pre-Hilbert $\mathcal{A}$-$\mathcal{B}$  bimodule. On completion, we obtain a genuine Hilbert $A$-$B$ bimodule which we denote by $E$. Moreover $E$ is an $A$-$B$ imprimitivity bimodule. For details, we refer the reader to \cite{MRW}. 

For a point $x \in X$, we denote the representation of $C^{*}(X \rtimes P)$ on $L^{2}(Q_x)$ induced at the point $x$ by  $\pi_{x}$ and the representation of $C^{*}(Y \rtimes \mathbb{R}^{d})$ on $L^{2}(\mathbb{R}^{d})$ induced at the point $x$ by $\widetilde{\pi}_{x}$. We claim that $Ind(\pi_x)=\widetilde{\pi}_x$.

Fix $x_0 \in X$. For $\chi \in C_{c}(Z)$, $\eta \in C_{c}(X \rtimes P)$ and $r \in \mathbb{R}^{d}$ , let 
\[
\widetilde{\chi \otimes \eta}(r)=\int \chi(x_0+r,s-r)\eta(x_0,s)1_{X}(x_0+s)ds.\]

\begin{enumerate}
\item[(1)] It is routine to  see that the map $C_{c}(Z) \otimes_{\mathcal{B}}C_{c}(X \rtimes P) \ni \chi \otimes \eta \to \widetilde{\chi \otimes \eta} \in C_{c}(\mathbb{R}^{d})$ is well-defined and extends to an isometry from    $E \otimes_{B} L^{2}(Q_{x_0})$ to the Hilbert space $L^{2}(\mathbb{R}^{d})$. 
\item[(2)] The set $\{\widetilde{\chi \otimes \eta}: \chi \in C_{c}(Z), \eta \in C_{c}(X \rtimes P)\}$ is total in $L^{2}(\mathbb{R}^{d})$.  Thus we can identify  $E \otimes_{B} L^{2}(Q_{x_0})$ with $L^{2}(\mathbb{R}^{d})$ via the unitary $E \otimes_{B} L^{2}(Q_{x_0}) \ni \chi \otimes \eta \to \widetilde{\chi \otimes \eta} \in L^{2}(\mathbb{R}^{d})$. Once  this identification is made,  a direct calculation shows that  $Ind(\pi_{x_0})=\widetilde{\pi}_{x_0}$. 
\end{enumerate}
Thus, in view of Theorem 6.23 of \cite{Rieffel} and Remark \ref{Rieffel commutant}, it suffices to prove Theorem \ref{main theorem of the appendix} for a transformation groupoid $Y \rtimes \mathbb{R}^{d}$. We do not claim any originality of what follows as it is well known. We include the details for completeness. 

Let us fix notation. Let $Y$ be a second countable, locally compact Hausdorff space on which $\mathbb{R}^{d}$ acts. Let $y \in Y$ be given. We denote the representation of $C^{*}(Y \rtimes \mathbb{R}^{d})$  induced at $y$ by $\pi_y$. Let $B(Y)$ be the algebra of bounded measurable functions on $Y$. For $f \in B(Y)$, let $M_{y}(f) \in B(L^{2}(\mathbb{R}^{d}))$ be defined by the equation
\[
M_{y}(f)\xi(t)=f(y+t)\xi(t)\]
for $\xi \in L^{2}(\mathbb{R}^{d})$. For $s \in \mathbb{R}^{d}$, let $L_{s}$ be the unitary on $L^{2}(\mathbb{R}^{d})$ defined by the equation
\[
L_{s}\xi(t)=\xi(t+s).\]
Then $(M_{y},L)$ is a covariant representation of the dynamical system $(C_{0}(Y),\mathbb{R}^{d})$.  If we identify $C^{*}(Y \rtimes \mathbb{R}^{d})$ with $C_{0}(Y) \rtimes \mathbb{R}^{d}$, the covariant representation that corresponds to the non-degenerate representation $\pi_y$ is $(M_y,L)$.

\begin{ppsn}
With the foregoing notation, we have the following.
\begin{enumerate}
\item[(1)] Let $y_1,y_2 \in Y$ be given. The representations $\pi_{y_1}$ and $\pi_{y_2}$ are non-disjoint if and only if there exists $s \in \mathbb{R}^{d}$ such that $y_1+s=y_2$. 
\item[(2)] For $y_0 \in Y$, the  commutant of the von Neumann algebra generated by the set $\{\pi_{y_0}(\xi): \xi \in C^{*}(Y \rtimes \mathbb{R}^{d})\}$ is the von Neumann algebra generated by $\{L_{s}: s \in H\}$ where $H$ is the stabiliser group of $y_0$, i.e. $H:=\{s \in \mathbb{R}^{d}: y_0+s=y_0\}$. 
\end{enumerate}
\end{ppsn}
\textit{Proof.} Suppose there exists $s \in \mathbb{R}^{d}$ such that $y_1+s=y_2$. Then $L_s$ intertwines $\pi_{y_1}$ and $\pi_{y_2}$. Conversely, suppose $\pi_{y_1}$ and $\pi_{y_2}$ are non-disjoint. Let $T$ be a non-zero intertwiner. Then $T$ intetwines $M_{y_1}$ and $M_{y_2}$. As a consequence, we have $TM_{y_1}(f)=M_{y_2}(f)T$ for every $f \in B(Y)$. For $i=1,2$, let $H_{i}$ be the stabiliser group of $y_i$. Note that the map $\mathbb{R}^{d}/H_{i} \ni s+H_{i} \to y_{i}+s \in Y$ is $1$-$1$ and continuous. We denote its image by $E_i$. Thanks to Theorem 3.3.2 of \cite{Arveson_invitation}, it follows that $E_{i}$ is a Borel set. Note that $M_{y_1}(1_{E_1})=1$. The equation $T=TM_{y_1}(1_{E_1})=M_{y_2}(1_{E_1})T$ implies that $M_{y_2}(1_{E_1}) \neq 0$. This implies that the orbit of $y_2$ meets the orbit of $y_1$. This proves $(1)$. 

It is clear that $\{L_{s}: s \in H\}$ lies in the commutant of $\pi_{y_0}(C_{0}(Y) \rtimes \mathbb{R}^{d})$. Conversely, suppose $T$ lies in the commutant of $\pi_{y_0}(C_{0}(Y) \rtimes \mathbb{R}^{d})$.  Then $T$ commutes with the algebra $\{M_{y_0}(f):f \in B(Y)\}$ and $\{L_{s}:s \in \mathbb{R}^{d}\}$.

Note that $\frac{\mathbb{R}^{d}}{H} \ni s+H \to y_0+s \in Y$ is continuous and $1$-$1$. Since $\mathbb{R}^{d}/H$ and $Y$ are Polish spaces, it follows from Theorem 3.3.2 of \cite{Arveson_invitation} that via the embedding $\mathbb{R}^{d}/H \ni s+H \to y_0+s \in Y$, we can identify the Borel space $\mathbb{R}^{d}/H$ with a subspace of $Y$. Thus bounded measurable functions on $\mathbb{R}^{d}/H$ can be considered as bounded measurable functions on $Y$. More precisely, suppose $f$ is a bounded measurable function on $\mathbb{R}^{d}/H$, then we consider $f$ as a function on $Y$  simply by declaring the values of $f$ on the complement of $\mathbb{R}^{d}/H$ to be zero. This way we   embedd  $C_{0}(\mathbb{R}^{d}/H)$ inside $B(Y)$. 

Hence we get a covariant representation $(M_{y_0},L)$ of $(C_{0}(\mathbb{R}^{d}/H), \mathbb{R}^{d})$. By Mackey's imprimitivity theorem, it follows that $C_{0}(\mathbb{R}^{d}/H) \rtimes \mathbb{R}^{d}$ is Morita equivalent to $C^{*}(H)$. Then the representation $M_{y_0} \rtimes L$ of $C_{0}(\mathbb{R}^{d}/H) \rtimes \mathbb{R}^{d}$ is $Ind(\rho)$ where $\rho$ is the regular representation of $H$ on $L^{2}(H)$. Note that $T$ lies in the commutant of $Ind(\rho)$. Since the group $H$ is abelian, it follows that the commutant of $\rho$ is the von Neumann algebra generated by $\{\rho(s):s \in H\}$. The statement now follows by appealing to Remark \ref{Rieffel commutant}. \hfill $\Box$

\nocite{Faraut}
\nocite{Arveson}
\nocite{Connes_Sur_la}
\nocite{Arveson_spectral}
\nocite{Rieffel}


\bibliography{references}
 \bibliographystyle{amsplain}

\noindent
 {\sc Anbu Arjunan}
(\texttt{aanbu@cmi.ac.in})\\
         {\footnotesize  Chennai Mathematical Institute, H1 Sipcot IT Park, \\
Siruseri, Padur, 603103, Tamilnadu, INDIA.}\\

\noindent
{\sc S. Sundar}
(\texttt{sundarsobers@gmail.com})\\
         {\footnotesize  Institute of Mathematical Sciences, CIT Campus, \\
Taramani, Chennai, 600113, Tamilnadu, INDIA.}\\

\end{document}